\documentclass{amsart}
\usepackage{amssymb}
\usepackage{multicol}
\usepackage{subfigure}
\usepackage{amscd}
\usepackage[usenames,dvipsnames]{pstricks}
\usepackage{epsfig}
\usepackage{pst-grad} 
\usepackage{pst-plot} 
\usepackage{verbatim}
\usepackage{enumerate}
\usepackage{verbatim}

\numberwithin{equation}{section}
\newtheorem{thm}{Theorem}[section]
\newtheorem{lemma}[thm]{Lemma}
\newtheorem{prop}[thm]{Proposition}

\def \re {\operatorname{Re}}
\def \im {\operatorname{Im}}

\def \mn {{\mathbb N}}

\def \eps {\varepsilon}   
\def \la {\lambda}

\def \ha{ {\frac{1}{2}}}

\def \p {\partial}

\def \rao#1 {\frac{\p}{\p #1} #1}

\newcommand{\sgn}{\operatorname{sgn}}

\def \mn {{\mathbb N}}

\def \eps {\epsilon}
\def \la {\lambda}

\def \p {\partial}

\def \rao#1 {\frac{\p}{\p #1} #1}

\usepackage{latexsym,eucal,amsfonts,amssymb,amsmath,graphicx}

\newtheorem{dfn}{Definition}

\newcommand{\beq}{\begin{equation}}
\newcommand{\eeq}{\end{equation}}

\def \mn {{\mathbb N}}

\def \eps {\varepsilon}
\def \la {\lambda}

\input epsf
\title{Resolvent Estimates on Asymptotically Hyperbolic Spaces}
\author{Raphael Hora}
\date{}
\begin{document}

\begin{abstract}
 We extend Vasy's results \cite{V2} on semiclassical high energy estimates for the meromorphic continuation of the resolvent for asymptotically hyperbolic manifolds to metrics that are not necessarily even. The Vasy's method \cite{Z} gives the meromorphic continuation of the resolvent and high energy estimates in strips, assuming that the geodesic flow is non-trapping, without having to construct a parametrix. We prove that the size of the strip to which the resolvent can be extended meromorphically is the same obtained by Guillarmou \cite{colin2}.
\end{abstract}

\maketitle
\section{Introduction}

 The particular class of manifolds studied here is modeled by the hyperbolic space 
 at infinity.  We equip the interior of a $C^\infty$ compact manifold $X$ of dimension  $(n+1)$ with boundary $\p X$ with a  $C^\infty$ metric $g$ such that $H=\rho^2 g$
is  non-degenerate up to $\partial X,$ where  $\rho$  defines  $\partial X$. 
As in \cite{mazzmel}, to make sure the sectional curvatures approach $-1$ at infinity, one also needs to assume that 
\begin{gather}
\displaystyle |d\rho|_{H}=1 \text{ at } \partial X. \label{asympt-curv}
\end{gather}
Such Riemannian manifolds $(X,g)$ are called asymptotically hyperbolic. 
It is shown in \cite{Gra} that if $g$ satisfies the conditions given above and $h_{0}=\rho^2 g|_{\p X}$  then there exist $\epsilon>0$ and a unique defining function of the boundary $x$ in $[0,\eps),$ such that  $\displaystyle X \sim [0,\epsilon)\times \partial X$  and 
\begin{equation}\label{metric2}
g=\frac{dx^2}{x^2}+\frac{h(x,y,dy)}{x^2},\quad h_{0}=h_{0}(0,y,dy).
\end{equation}
One can think of $h(x,y,dy)=h(x)$ as a $C^{\infty}$ one-parameter family of metrics on the boundary $\partial X$. We fix this decomposition, and from now on, $x\in C^{\infty}(X)$ is as in \eqref{metric2}.   Since $h(x,y,dy)$ is $C^\infty,$ it has a Taylor series expansion at $x=0:$ 
\begin{gather*}
h(x,y,dy)\sim h_0(y,dy)+ \sum_{j=1}^N h_{j}(y,dy)x^{j} + O(x^{N+1}), \;\ N\in \mn,
\end{gather*}
and following Guillarmou \cite{colin2} we say that $g$ is even to order $O(x^{2k+1})$ if  $h_j(y,dy)=0$ for all odd values of $j=2m+1,$ with $m<k.$  Throughout this paper we will assume that $g$ is even modulo $O(x^{2k+1}),$ $k\geq 2.$ 


The spectrum of $\Delta_{g}$ was studied by Mazzeo and Melrose in \cite{mazzeo1,mazzeo2,mazzmel} and more recently by Bouclet \cite{bouclet}. It consists of a finite point spectrum $\sigma_{pp}(\Delta_{g})$, which is the set of $L^{2}(X)$ eigenvalues, and an absolutely continuous spectrum $\sigma_{ac}(\Delta_{g})$ satisfying
\begin{equation}\label{l2decomp}
\sigma_{ac}(\Delta_{g})=[n^2/4,\infty),\quad \sigma_{pp}(\Delta_{g})\subset(0,n^2/4).
\end{equation}
It follows from the spectral theorem that if $\text{Im }\lambda>>0$, the resolvent for $\Delta_{g}$, denoted by
\begin{equation}
R(\lambda)=\left(\Delta_{g}-\frac{n^2}{4}-\lambda^2\right)^{-1}
\end{equation}
is a bounded operator in $L^{2}(X).$ Mazzeo and Melrose also proved in \cite{mazzmel} that, as an operator $R(\lambda):\mathcal{\dot{S}}(X)\rightarrow \mathcal{\dot{S}}'(X),$ the resolvent has a meromorphic extention to $\mathbb{C}\setminus(-\frac{i}{2}\mathbb{N})$, where $\displaystyle \mathcal{\dot{S}}(X)$ is the space of Schwartz functions on $X$ which vanish with all derivatives at $\partial X$, and $\displaystyle \mathcal{\dot{S}}'(X)$ is its dual space of tempered distributions. Guillarmou showed in \cite{colin2} that  if the metric is even to order $O(x^{2k+1})$ then
the resolvent has a finite meromorphic extension to a half-space $\{z\in \mathbb{C}:\ \im z>  -1/2-k\}$, as an operator acting on weighted $L^2$ spaces. Conversely, if  $k\geq 2$ and $R(\la)$ extends to $\{z\in \mathbb{C}\;\ \im z>- 1/2-k\}$, then the metric has to be even to order $O(x^{2k-1}).$ In particular  $R(\la)$ has a meromorphic extension to $\mathbb{C}$ if and only if $g$ is even to all orders.  Guillarmou  also showed that if $g$ is not even, then generically $R(\lambda)$  has essential singularities on $\frac{-i}{2}\mathbb{N}.$

When $(X,g)$ is a perturbation of the Poincar\'e metric on hyperbolic space supported in a small neighborhood of $\p X,$ but the metrics is not necessarily even,  Melrose, S\'a Barreto and Vasy \cite{MVS2} proved  high energy semiclassical estimates for the resolvent operator on strips, and away from the imaginary axis. Chen and Hassel \cite{CheHa} and Wang \cite{Wang} generalized the results of \cite{MVS2} for non-trapping AHM. The case of conformally compact manifolds with variable negative sectional curvature at infinity has also been studied. Borthwick \cite{Bo} constructed a parametrix to obtain the meromorphic continuation of the resolvent and S\'a Barreto and Wang \cite{SaWa} constructed a semiclassical parametrix of the resolvent and used it to establish high energy resolvent estimates. Estimates just on the real axis were obtained by Cardoso and Vodev for a more general class of metrics and Guillarmou obtained in \cite{colin1} high energy estimates on an exponentially small neighborhood of the real axis for general metrics, and on a strip for non-trapping metrics of constant curvature near infinity. 

We will denote by $H^{s}(X)$ the standard Sobolev space of order $s$, and by $\overline{H}^{s}(X)$ corresponding to extensions of $H^{s}$ functions across the boundary. We denote by $H^{s}_{\hbar}(X)$ the semiclassical Sobolev space. We recall that for $h$ bounded away from $0$, if $s\in\mathbb{N},$ the $H^{s}_{\hbar}$ norm is locally given by 
$$ \|u\|^{2}_{H^{s}_{\hbar}}=\sum_{|\alpha|\leq s}\|(hD_{j})^{s}u\|_{L^2}^2.$$
More generally $\displaystyle \|u\|^{2}_{H^{s}_{\hbar}}=\|(1+|\xi|^2)^{s}\mathcal{F}_{\hbar}u\|_{L^2}^2,$ where $\mathcal{F}_{\hbar}$ is the semiclassical Fourier transform. In $\mathbb{R}^n$, $\mathcal{F}_{\hbar}$ is given by
$$\mathcal{F}_{\hbar}u(\xi)=\int_{\mathbb{R}^n}e^{-ix\cdot \xi/h}u(x)dx.$$

We extend Theorem 5.1 of \cite{V2} to the class of metrics $g$ which are even of order $O(x^{2k+1}):$
\begin{thm}\label{mainthm1}\label{mainthm}
Let $(X,g)$ be an asymptotically hyperbolic manifold, and $k\in \mathbb{N}\cup\{\infty\}.$ If $g$ is even modulo $O(x^{2k+1}), k\geq2,$ the modified resolvent 
$$R(\lambda)=\left(\Delta_g-\frac{n^2}{4}-\lambda^2\right)^{-1}:\mathcal{\dot{S}}\longrightarrow \mathcal{\dot{S}}'$$
defined for $\text{Im }\lambda>>0$, extends to a finite-meromorphic family from $\text{Im }\lambda>>0$ to $\text{Im }\lambda>-1/2-k.$ Moreover if the metric $g$ is assumed to be non-trapping, there exists $\epsilon_0>0$ such that in $1/2-s<-C<\text{Im }\lambda<\epsilon_0|\text{Re }\lambda|,$ $|\text{Re }\lambda|>C$ non-trapping semiclassical estimates hold for $-k<s< k+1$
\begin{equation}\label{nontrapest3}
\|x^{-n/2}e^{i\lambda\phi}R(\lambda)f\|_{H_{|\lambda|^{-1}}^s}\leq C'|\lambda|^{-1}\|x^{-(n+4)/2}e^{i\lambda\phi}f\|_{H_{|\lambda|^{-1}}^{s-1}},
\end{equation}
where $\displaystyle e^{\phi}=x(1+x^2)^{-1/4}$ near $\{x=0\}.$ In particular, for $s=1,$ $\im \lambda>-C>-1/2,$ $|\re \lambda|>C$ (\ref{nontrapest3}) is
$$\|x^{-n/2}e^{i\lambda\phi}R(\lambda)f\|_{L^2}+|\lambda|^{-2}\|d(x^{-n/2}e^{i\lambda\phi}R(\lambda)f)\|_{L^2}\leq C'' |\lambda|^{-2}\|x^{-(n+4)/2}e^{i\lambda\phi}f\|_{L^2}.$$ 
\end{thm}

\begin{figure}[h!]
\begin{center}
\includegraphics[scale=0.40]{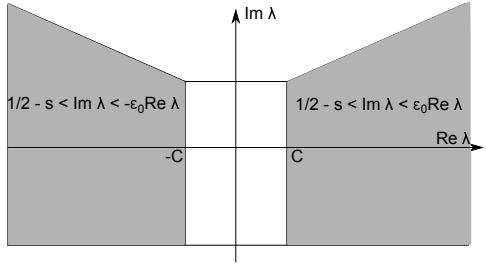}
\caption{The shaded region shows where the resolvent estimate is valid. Notice that the physical plane here is $\im \la>0$ and  we analyze the continuation of $R(\la)$ to $\im \la<0.$ }
\end{center} 
\end{figure}

As noted by Vasy \cite{V2}, if one assumes that $f$ has compact support in $X$, the $s-1$ norm on $f$ can be replaced by the $s-2$ norm, since the operator is semiclassically elliptic. Also, the resolvent is semiclassically outgoing in the sense of Datchev and Vasy \cite{DaVa1}. Then one obtains resolvent estimates similar to (\ref{nontrapest3}) if one assumes mild trapping such as normally hyperbolic trapped sets. 

 The first step of the  proof is to change the $C^\infty$ structure of $X$ by making $\mu=x^2,$ which means  that only even functions of $x$ are smooth in $\mu$. If the metric $g$ is even,  the  Laplacian $\Delta_{g}$ remains a $C^\infty$ operator in the new $C^\infty$ structure. In the case of even metrics of order $O(x^{2k+1}),k\geq2,$ $\Delta_g$ will have limited smoothness and the novelty of this paper is to adapt the methods of \cite{V2} to the case of non-smooth coefficients by  working with non-smooth pseudodifferential operators. 

This gives another perspective of the results of  Mazzeo and Melrose \cite{mazzmel}. We also give a new explanation and recover Guillarmou's results \cite{colin2} and also give high energy estimates. The estimate \eqref{nontrapest3} is stronger than the one obtained by Melrose, S\'a Barreto and Vasy in \cite{MVS2}. As mentioned by Zworski \cite{Z}, Vasy's method gives a rather direct way to prove the meromorphic continuation of the resolvent of the Laplacian on asymptotically hyperbolic manifolds and also gives semiclassical resolvent bounds by relating the resolvent to the inverse of a family of Fredholm differential operators. 

This paper was part of my Ph.D thesis, and I would like to thank my  advisor Ant\^onio S\'a Barreto for all the help. I am also grateful to the NSF for their support under grant DMS 0901334.

\section{Outline of the Proof}

We follow the methods of  \cite{V2}. As mentioned above, the first step is to change the $C^{\infty}$ structure of $X$, by setting $\mu=x^2,$ that is, $x=\sqrt{|\mu|}.$ We denote the manifold with the new $C^\infty$ structure by $X'.$ We consider a conjugated operator $\displaystyle P_{\lambda}=e^{i\lambda\varphi}(\Delta_{g}-n^2/4-\lambda^2)e^{-i\lambda\varphi}$, where $\displaystyle e^{\varphi}=\mu^{1/2}(1+\mu)^{-1/4}$ in a neighborhood of $\partial X'=\{\mu=0\}$.  It turns out that $P_{\lambda}$ can be extended to the exterior of the manifold $X$. We show that the conjugated operator $P_{\lambda}$ is semiclassically elliptic in the interior of $X'$ ($\mu<0$) if $\varphi$ satisfies $4\mu^2(\partial_{\mu}\varphi)^2<1.$  $\displaystyle P_{\lambda}$ is not smooth at $\{\mu=0\},$ but since $k\geq 2,$ we are able to also extend it across the boundary, i.e.  for $\mu<0$. 

Since $k\geq 2,$ the derivatives of the principal symbol of $\displaystyle P_{\lambda}$ are Lipschitz, and hence the dynamics of the hamiltonian flow of the principal symbol of the conjugated operator is not affected by the order of evenness of the metric.  We work in the pre-compact manifold without boundary $X_{-\delta_0}=\left((-\delta_{0},0)_{\mu}\times \partial X\right)\cup X'.$ $P_{\lambda}$ is not semicassically elliptic if $\mu\leq 0$. We then show that the characteristic set of the conjugated operator splits into two disjoint parts. We say that the characteristic set has a negative and a positive part, depending on the sign of $\xi$ (the dual variable of $\mu$). The negative part is a sink and the positive one is a source, in the sense that all nearby bicharacteristics converge to the positive and negative parts respectively as the time variable goes to $\mp\infty$. The null bicharacteristics of $P_\lambda$ are radial at the conormal bundle $N^{*}\partial X'\setminus o$ of $\partial X'=\{\mu=0\}.$ Vasy \cite{V2} adapted Melrose's Theorem \cite{mesc} on propagation of singularities at radial points to prove an estimate for $P_\lambda$ near $\{\mu=0\}.$ We then extend these results to our setting by working with a calculus of pseudodifferential operators with non-smooth coefficients. 

To deal with the region where $\mu<0$, we introduce a complex absorption operator $Q_{\lambda}$ so that $P_{\lambda}- iQ_{\lambda}$ is elliptic near $\mu=-\delta_0$ and $Q_{\lambda}$ is supported in $\mu<0$. As it is shown in Appendix E.5 of the book of Dyatlov and Zworski \cite{DyZw}, this can be arranged if $q$ has the correct sign along the sink and the source. We show that,  for $s>m>1/2-\text{Im }\lambda,$ $-k<m<s<k+1,$ if $u\in H^{m}(X_{-\delta_0})$ and $(P_{\lambda}-iQ_{\lambda})u\in H^{s-1}(X_{-\delta_0})$, then $u\in H^{s}(X_{-\delta_0})$. Instead of considering an complex absorption operator, Zworski \cite{Z} showed that the operator $P_\lambda$ is hyperbolic when $\mu<0$, then one can use standard hyperbolic estimates. The regularity results for $P_\lambda$, namely propagation of singularities, will give a bound of $\displaystyle \|u\|_{H^s}$ in terms of $\|(P_\lambda-iQ_\lambda)u\|_{H^{s-1}}$ (instead of $\|(P_\lambda-iQ_\lambda)u\|_{H^{s-2}}$). 

It then makes sense to consider the map 
$$i_{s}:E_{s}=\{u\in H^{m}:(P_{\lambda}-iQ_{\lambda})u\in H^{s-1}(X_{-\delta_0})\}\rightarrow H^{s}(X_{-\delta_0}).$$ 
We show that $E_{s}$ is a Banach space and as an application of the closed graph theorem, $i_{s}:E_{s}\rightarrow H^{s}(X_{-\delta_0})$ is a continuous map. Since we are now working on a pre-compact manifold without boundary, the map $i_{s}$ is compact, which implies that $\text{Ker}(P_{\lambda}-iQ_{\lambda})$ is finite-dimensional. Similarly, we prove that $\text{Ker}(P_{\lambda}^{*}+iQ_{\lambda}^{*})$ also has finite dimension. We then follow the ideas of H\"ormander given in the proof of Theorem 26.1.7 of \cite{Ho} to show that for any $f$ in the annihilator of $\text{Ker}(P_{\lambda}^{*}+iQ_{\lambda}^{*})$, there exists $u$ such that
$$(P_{\lambda}-iQ_{\lambda})u=f, \quad \text{in } X.$$
This $H^{s}-$solvability combined with the fact that the kernel and cokernel of $P_{\la}-iQ_{\la}$ are finite-dimensional implies that $P_{\la}-iQ_{\la}$ is an analytic family of Fredholm operators. We then apply analytic Fredholm theory to show that this operator has a bounded inverse for $s>1/2-\text{Im }\lambda$, $-k<s\leq k+1.$ In \cite{V2}, since the metric is even there is no constraint in $s$ with respect to $k$ and Vasy showed meromorphic continuation of $R(\lambda)$ to $\mathbb{C}$, but as proved by Guillarmou in \cite{colin2}, if the metric is even of order $O(x^{2k+1})$, $R(\lambda)$ continues meromorphically only to a strip which depends on $k.$ 

Finally, we show the same meromorphic continuation for the respective semiclassical operators (We call them $P_{\lambda,\hbar}$ and $Q_{\lambda,\hbar}$, where $h$ is the semiclassical parameter). As in \cite{V2}, the boundedness of the operator $(P_{\lambda,\hbar}-iQ_{\lambda,\hbar})^{-1}$ in $H^{s}$, for $s>1/2-\text{Im }\lambda$, $-k<s\leq k+1,$ implies the estimate (\ref{nontrapest3}). We finally relate the operator $(P_{\lambda,\hbar}-iQ_{\lambda,\hbar})^{-1}$ with the resolvent $R(\lambda)=(\Delta_{g}-n^2/4-\lambda^2)^{-1}$ by the properties of the support of the complex absorption operator, namely its support is in $\{\mu<0\}$.

\section{The extended operator}

We will state and repeat some of Vasy's ideas for the covenience of the reader. Suppose that $g$ is an asymptotically hyperbolic metric on $X$ with $\text{dim }X=n+1.$ It follows from \eqref{metric2} that, in a collar neighborhood of the boundary, the Laplace-Beltrami operator with respect to $g$ is given by

\begin{equation}
\Delta_{g}=(xD_{x})^2+i(n+x^2\gamma)(xD_{x})+x^2\Delta_{h},
\end{equation}
with $\displaystyle\gamma=-\frac{1}{2x}D_{x}\log|h|$ even modulo $O(x^{2k-1})$, and $\Delta_{h}$ the $x-$dependent family of Laplacians of $h$ on $\partial X.$

If we introduce $\mu=x^2$ as the new boundary defining function, this changes the $C^{\infty}$ structure of $X$ and only even fuctions of $x$ are smooth. We have then 
\begin{equation}
\Delta_{g}=4(\mu D_{\mu})^2+2i(n+\mu\gamma)(\mu D_{\mu})+\mu \Delta_{h}.
\end{equation}

Now if we conjugate by $\mu^{-i\lambda/2+(n+2)/4}$ and multiply by $\mu^{-1/2}$ from both sides, we obtain
\begin{equation}\label{opvasy}
\begin{gathered}
\mu^{-1/2}\mu^{i\lambda/2-(n+2)/4}\left(\Delta_{g}-\frac{n^2}{4}-\lambda^2\right)\mu^{-i\lambda/2+(n+2)/4}\mu^{-1/2}\\
=4\mu D_{\mu}^2-4\lambda D_{\mu}+\Delta_{h}-4iD_{\mu}+2i\gamma\left(\mu D_{\mu}-\frac{\lambda}{2}-i\frac{n}{4}\right).
\end{gathered}
\end{equation}
Since $g$ is assumed to be even modulo $O(x^{2k+1})$, some coefficients of terms order lower than two of the operator (\ref{opvasy}) above are only smooth modulo $O(\mu^{k-1/2})$, so the operator in (\ref{opvasy}) is in $\displaystyle C^{k-1}_{\mu}\text{Diff}^2(X')$ (that is, a differential operator of order two with $C^{k-1}$ coefficients with respect to the $\mu$ variable), but the other coefficients are $C^{k}$ with respect to $\mu$. $P_\lambda$ continues across the boundary, by extending $h$ and $\gamma$ in an arbitrary $C_{\mu}^{k-1}$ manner. 

As noted by Vasy in \cite{V2}, in order to archieve an operator which is semiclassically elliptic when $\lambda$ is away from the real axis we have to make one more conjugation. We then conjugate by $(1+\mu)^{i\lambda/4}$ and obtain
\begin{equation}\label{psig}
P_{\lambda}=4\mu D_{\mu}^2-4(1+a_{2})\lambda D_{\mu}-\lambda^2+\Delta_{h}-4i D_{\mu}+b_{1}\mu D_{\mu}+b_2 \lambda+c_1
\end{equation}
with 
\begin{gather*}
a_{2}=-\frac{\mu}{2(1+\mu)}, \quad b_{1}=2i\gamma,\quad \\
 b_{2}=\frac{1}{1+\mu}\left[\left(\frac{\lambda}{4}+i\right)\frac{\mu}{1+\mu}-(1+i)-\frac{i\gamma(2+\mu)}{2}\right],\quad 
c_1=\frac{(n-1)}{2}\gamma.
\end{gather*}

The Taylor expansion of $h$ in $\mu$ is
$$h=h_0+\mu h_{2}+\mu^{2}h_{4}+\ldots \mu^k h_{2k} + |\mu|^{k+\ha} h_{2k+1}+ \ldots \;\ \mu\geq 0,$$
so  $h$ can be continued  to $X_{-\delta_0}=(-\delta_0,0)\times \partial X\cup X',$ and hence we can continue the operator $P_{\lambda}$ defined by \eqref{psig} for $\mu<-\delta_0<0.$  It also follows from \eqref{psig} that  the extension of $P_{\lambda}$ is non-degenerate.

Writing covectors as $\xi d\mu+\eta dy,$ the standard classical principal symbol of $P_{\lambda}$ is given by
\begin{equation}\label{princesymb}
p=\sigma_{2}(P_{\lambda})=4\mu\xi^2+|\eta|_{\mu,y}^2, \text{ where } |\eta|_{\mu,y}^2=\sum_{1\leq j,k\leq n}h^{jk}\eta_{j}\eta_{k},
\end{equation}
which is real, independent of $\lambda$. 

Recall that if an operator is formally self-adjoint relative to the metric density $|dg|$, then its conjugate by $f$ ($f$ real function) is formally self-adjoint with respect to $f^{2}|dg|$. Since $\displaystyle \Delta_{g}-\lambda^2-\frac{n^2}{4}$ is formally self-adjoint relative to $|dg|$, then $P_{\lambda}$ is formally self-adjoint with respect to $\mu^{(n+2)/2}|dg|=\frac{1}{2}|dh||d\mu|,$ as $\displaystyle x^{-n-1}dx=\frac{1}{2}\mu^{-(n+2)/2}d\mu$. Therefore $\displaystyle \mu^{(n+2)/2}|dg|$ extends to a density to $X_{-\delta_0}$.

We can also conclude that $\displaystyle \sigma_{1}\left(\frac{1}{2i}(P_{\lambda}-P_{\lambda}^{*})\right)\left|_{\mu\geq0}\right.$ vanishes if $\lambda\in\mathbb{R}.$ Indeed, we have that 
\begin{equation}\label{difop1}
\sigma_{1}(P_{\lambda}-P_{\lambda}^{*})|_{\mu=0}=-8i(\text{Im }\lambda)\xi.
\end{equation}

Finally we need to check that $\mu$ can be appropriately chosen in the interior of $X$, that is $\mu>0$. As noted by Vasy in \cite{V2}, this only matters for semiclassical purposes, and for $z$ non-real, the operator $\displaystyle e^{i\lambda\varphi}(\Delta_{g}-n^2/4-\lambda^2)e^{-i\lambda\varphi}$ is semiclassically elliptic in the interior of $X$ if we extend $\varphi$ so that $\displaystyle 4\mu^2(\partial_{\mu}\varphi)^2<1$ and $e^{\phi}=\mu^{1/2}(1+\mu)^{-1/4}$ near $\{\mu=0\}$. Then, after conjugation by $e^{-i\lambda\varphi}$,
\begin{equation}
P_{h,z}=e^{iz\varphi/h}\mu^{-(n+2)/4-1/2}\left(h^2\Delta_{h}-z^2-\frac{n^2}{4}\right)\mu^{(n+2)/4-1/2}e^{-iz\varphi/h}
\end{equation}
is semiclassically elliptic in $\mu>0,$ as we wanted.

\section{Pseudodifferential operators with nonsmooth symbols}
Since the metric $g$ is even of order $O(x^{2k+1})$, the coefficients of $P_{\lambda}$ are $C^{k-1}$ with respect to $\mu$. We then have to consider a calculus of pseudoffirential operators with non-smooth symbols. We recall some results from \cite{BR1,BR2,Ma1,Ta1,Ta3}  for convenience. The calculus for these operators with singular symbols is analogous to the classical calculus of pseudodifferential operators except that, because of the limited  regularity, expansions have only finitely many terms. It is clear that if the coefficients are non-smooth, then they themselves will introduce singularities into the solutions of $P_{\lambda}u=f$ so that only a limited amount of smoothness in the initial data can be propagated to later times. All the definitions given below are in euclidean spaces, but as Marschall showed in \cite{Ma1}, the classes defined in this section are invariant under coordinate transformation and one can easily transfer them to a manifold $X$ by using coordinate charts.
\begin{dfn}
 Define $S^{m}(r)$ to consist of symbols $a:\mathbb{R}_{\mu}\times\mathbb{R}^n_{y}\times\mathbb{R}_{\xi}\times\mathbb{R}^n_{\eta}\rightarrow\mathbb{C}$ such that for each $\displaystyle \alpha\in\mathbb{N}$ and $\gamma,\beta\in\mathbb{N}^{n},$ 
\begin{equation}
\begin{gathered}
|\partial_{y}^{\gamma}\partial^{\beta}_{\eta}\partial^{\alpha}_{\xi}a(\mu,y,\xi,\eta)|\leq C_{\alpha,\beta,\gamma}(1+|\xi|+|\eta|)^{m-\alpha-|\beta|},\\
\|\partial^{\gamma}_{y}\partial^{\beta}_{\eta}\partial^{\alpha}_{\xi}a(\cdot,y,\xi,\eta)\|_{H^{r}}\leq C_{\alpha,\beta,\gamma}(1+|\xi|+|\eta|)^{m-\alpha-|\beta|},
\end{gathered}
\end{equation}
where $\|\cdot\|_{H^{r}}$ is the Hardy-Sobolev norm with respect to the first coordinate $\mu$, that is, given the Bessel potential of order $k$
$$J^{k}f(\mu)=\frac{1}{2\pi}\int e^{i\mu\xi}(1+|\xi|^2)^{k/2}\hat{f}(\xi)d\xi,$$
and a Schwartz function $\phi$ such that $\displaystyle\int \phi d\mu=1$, let $\phi_{t}(x)=t^{-1}\phi(x/t),$ then
$$\|f\|_{H^r}=\|J^{r}f\|_{h^2}=\left\|\sup_{0<t<1}|\phi_t\ast J^{r}f|\right\|_{L^2}.$$
\end{dfn}
We write just $a\in S^{l}$ if $a\in S^{l}(\infty),$ that is $a\in S^{l}(r)$ for all $r$. We now define the class of standard symbols, which includes the ones that we are considered in this work.

\begin{dfn}\label{stps}
For $k\geq 0$ an integer, $\displaystyle S_{st}^{m+k}(r+k)$ consists of the symbols $a(\mu,y,\xi,\eta)$ of the form
$$a(\mu,y,\xi,\eta)=a_{m+k}(\mu,y,\xi,\eta)+a_{m+k-1}(\mu,y,\xi,\eta)+\cdot\cdot\cdot+a_{m}(\mu,y,\xi,\eta),$$
where $a_{m}(\mu,y,\xi,\eta)\in S^{m}(r)$ and, for $j=1,...,k,$ 
$$a_{m+j}(\mu,y,\xi,\eta)=\sum_{l}a_{j,l}(\mu,y)p_{m+j,l}(\mu,y,\xi,\eta)$$
with $a_{j,l}\in H_{\mu}^{r+j}$ and $p_{m+j,l}\in S^{m+j}.$
\end{dfn}

In a collar neighborhood of the boundary, the full symbol of our operator $P_{\lambda}$ is of the following form:

\begin{equation}
p_{\text{full}}=4\mu\xi^2-4(1+a_2)\lambda\xi-\lambda^2+\sum_{i,j}h^{ij}\eta_i\eta_j+\sum_{i,j}\left[(\partial_{y_i}log|h|)h^{ij}+\partial_{y_i}h^{ij}\right]\eta_j-4i\xi+b_1\mu\xi+b_2\lambda+c_1,
\end{equation}
where $a_2\in C^{\infty}_{\mu},$ $b_{1},b_{2},c_1\in C^{k-1}_{\mu}$, $\displaystyle \left[(\partial_{y_i}|h|)h^{ij}+\partial_{y_i}h^{ij}\right]\in C^{k-1}_{\mu}$ and $h^{ij}\in C^{k}_{\mu}.$  Then $p_{\text{full}}\in S_{st}^{2}\left(k\right),$ since the principal symbol $\displaystyle p=4\mu\xi+\sum_{i,j}h^{ij}\eta_i\eta_j$ has coefficients in $H^{k}$ (with respect to $\mu$). 
 
Given $a\in S^{m}(r)$, we define a pseudodifferential operator $A\in \Psi^{m}(r)$ by
\begin{equation}
\begin{gathered}
(A u)(z)=(Op(a)u)(z)=(2\pi)^{-n}\int e^{iz\cdot \tau}a(z,\tau)\hat{u}(\tau)d\tau,\\
\text{where } z=(\mu,y),\tau=(\xi,\eta).
\end{gathered}
\end{equation}

We will also use the notation $Op(a)=a(z,D)$ or $Op(a)=a(\mu,y,D_{\mu},D_{y})$, where $D$ represents the derivative with respect to $z=(\mu,y)$, that is $\xi\sim D_{\mu}$ and $\eta\sim D_{y}$. 

For a standard symbol $a\in S_{st}^{m+k}(k),$ $A=Op(a)\in \Psi^{m+k}(k)$, the principal symbol $\sigma_{m+k}(A)$ is\\ $a_{m+k}(\mu,y,\xi,\eta)$ as in Definition \ref{stps}. In particular if $A$ is a differential operator with rough coefficients, that is, non-smooth coefficients, say 
$$A=\sum_{|\alpha|\leq m}a_{\alpha}(z)D_{z}^{\alpha},$$
where $\displaystyle a_{\alpha}\in C^{k},$ then the principal symbol of $A$ is
$$\sigma_{m}(A)=\sum_{|\alpha|=m}a_{\alpha}(z)\zeta^{\alpha}.$$

When composing two pseudodifferential operators with singular symbols, say $a(\mu,y,D_{\mu},D_{y})$ and\\
$b(\mu,y,D_\mu,D_y)$, the usual complete asymptotic expansion 
$$a\circ b \sim \sum_{\alpha=(\alpha_1,\alpha')}(i^{-|\alpha|}/\alpha!)\partial_{\xi}^{\alpha_1}\partial_{\eta}^{\alpha'}a\partial_{y}^{\alpha'}\partial_{\mu}^{\alpha_1}b,$$
does not make sense, since  if  $\alpha_1 $ large the derivatives $\displaystyle\partial_{\mu}^{\alpha_1}$ are not defined. Actually, one can see that if $a\in S^m(r),$ and $\displaystyle\alpha_1<r$, then $\displaystyle D_{\mu}^{\alpha_1}a\in S^{m}(r-\alpha_1).$ Therefore we break off the expansion at step $k$, and set
\begin{equation}
\begin{gathered}
r_{k}(\mu,y,\xi,\eta)=\sigma_{A\circ B}(\mu,y,\xi,\eta)
-\sum_{\alpha=(\alpha_1,\alpha'),|\alpha|\leq k}(i^{-|\alpha|}/\alpha!)\partial_{\xi}^{\alpha_1}\partial_{\eta}^{\alpha'}a(\mu,y,\xi,\eta)\partial_{y}^{\alpha'}\partial_{\mu}^{\alpha_1}b(\mu,y,\xi,\eta).
\end{gathered}
\end{equation}
To be more precise, we have the following result that was proved in \cite{Ma1} (Corollary 3.4). 
\begin{lemma}\label{compoab}\label{abfou}\label{corcomp}\label{commubd}

Let $a\in S^{m_1}(r)$ and $b\in S^{m_2}(r)$ and let $\tau+1/2<r,$ $\tau\leq 1.$ Then $ab\in S^{m_1+m_2}(r)$ and for each $s$ such that
$$\tau-r+\max\{0,-m_1\}<s\leq r-\max\{0,m_1\}$$
the commutator $\displaystyle Op(a) Op(b)-Op(ab):H_{loc}^{s+m_1+m_2-\tau}\rightarrow H_{loc}^{s}$ is bounded.
\end{lemma}

We use the lemma above when $r=k\geq 2$, so we can take $\tau=1. $ This lemma can be generalized to the following result (similar to Lemma 1.5 of \cite{BR2}).

\begin{lemma}\label{comporem}
Let $s>1/2$, $j\geq0$ be an integer, and assume that $b\in S^{m}(l+j)$ and $a\in S^{j}(l)$. Then $Op(a)Op(b)$ is an operator with symbol $a\circ b\in S^{m+j}(l)$ satisfying 
$$r_{j}(\mu,y,\xi,\eta)=a\circ b(\mu,y,\xi,\eta)-\sum_{\alpha=(\alpha_1,\alpha'),|\alpha|< j}(i^{-|\alpha|}/\alpha!)\partial_{\xi}^{\alpha_1}\partial_{\eta}^{\alpha'}a(\mu,y,\xi,\eta)\partial_{y}^{\alpha'}\partial_{\mu}^{\alpha_1}b(\mu,y,\xi,\eta)\in S^{m}(l)$$ 

\end{lemma}

The following lemma (Corollary 3.6 in \cite{Ma1}) follows from an argument used in the standard calculus of smooth pseudodifferential operators and lemma \ref{abfou} above.

\begin{lemma}
Let  $a\in S^{m}(r)$ and $0\leq\tau\leq 1$ be such that $\tau+1/2<r$. Then for $\tau-r+\max\{0,-m\}<s<r-\max\{0,m\}$ (if $m>0$, then we can include $s=r-m$), the commutator $\displaystyle Op(a)^{*}-Op(\overline{a}):H^{s+m-\tau}\rightarrow H^{s}$ is bounded. 

\end{lemma}



We also need operator estimates on Sobolev spaces.  
\begin{prop}\label{TayOp}
Let $a\in S^{m}(r)$. Then, if $s$ satisfies $\displaystyle -r<s\leq r$, the operator $\displaystyle A:H_{loc}^{s+m}\rightarrow H_{loc}^{s}$ is bounded. 
\end{prop}

The following propagation result is proved in chapter II of \cite{tay} and Theorem 4.6 of \cite{Ma1}.   

\begin{prop}\label{ellipticregu}
Assume that $r>0$, $p(\mu,y,\xi,\eta)\in S^{m}(r)$ is elliptic, $-r<\delta<s<r$ and that 
$$u\in H_{loc}^{m+\delta},\quad p(\mu,y,D_{\mu},D_{y})u\in H_{loc}^{s}.$$
Then 
$$u\in H_{loc}^{s+m}.$$
\end{prop}

We can also state the Proposition above microlocally as: if $p(\mu,y,\xi,\eta)\in S^{m}(r)$ is elliptic near $\nu$, then if there exists a $\delta$ such that $-r<\delta<s<r$ and $\nu\not\in WF^{\delta+m}(u)$, we have
$$\nu\not\in WF^{s}(p(\mu,y,D_{\mu},D_{y})u)\Rightarrow \nu\not\in WF^{s+m}(u).$$

One important estimate is the sharp G\r{a}rding inequality (proposition 2.4.A in \cite{Ta4}).

\begin{prop}\label{Garding}
Let $p(\mu,y,\xi,\eta)\in S^{m}(l)$ and $p(\mu,y,\xi,\eta)\geq -C_0.$ Then for all $u\in C^{\infty}_{0},$ 
\begin{equation}
\text{Re }\langle p(\mu,y,D_\mu,D_y)u,u\rangle\geq -C_1\|u\|_{L^2}^2,
\end{equation}
provided $l>0$ and $\displaystyle m\leq \frac{2l}{2+l}.$
\end{prop}

Note that this proposition can be written as: let $p(\mu,y,\xi,\eta)\in S^{m}(l),$ then $$\tilde{p}(\mu,y,\xi,\eta)=J^{-r}p(\mu,y,\xi,\eta)J^{-r}\in S^{m-2r}(l+r),$$ that is $p(\mu,y,\xi,\eta)=J^r\tilde{p}(\mu,y,\xi,\eta)J^r\in S^m(l).$ So if $p(\mu,y,\xi,\eta)\geq -C_0,$ then for all $u\in C^{\infty}_{0}$,
\begin{equation}
Re\langle p(\mu,y,D_\mu,D_y)u,u\rangle\geq -C_1\|u\|_{H^r}^2,
\end{equation}
provided $l+r>0$ and $\displaystyle m-2r\leq \frac{2l+2r}{2+l+r}.$



The semiclassical pseudodifferential algebra is analogous. We say that $A\in \Psi_{\hbar}^{m}(k)$ if it is given by
\begin{equation}
(A_{h}u)(z)=(2\pi h)^{-n}\int_{\mathbb{R}^n}e^{i(z-z')\cdot\tau/h}a(z,\tau,h)u(z')d\tau dz',
\end{equation}
where $u\in \mathcal{S}(\mathbb{R}^n)$ and $a\in C^{\infty}\left([0,1)_{h};S^{m}(k)\right)$, $z=(\mu,y).$ We also write $a_{\hbar}(z,\tau)=a(z,\tau,h),$ and $a_{\hbar}\in S_{\hbar}^{m}(k)$ instead of $a\in C^{\infty}\left([0,1)_{h};S^{m}(k)\right)$. $A_{\hbar}$ is called the quantization of the symbol $a_{\hbar}$. Then a semiclassical differential operator (with an extra variable $\lambda$ as in our case) is written as
$$A_{h,\lambda}=\sum_{|\alpha|\leq m}a_{\alpha}(z,\lambda;h)(hD_{z})^{\alpha},$$
with semiclassical principal symbol
$$\sigma_{\hbar,m}(A_{h,\lambda})=\sum_{|\alpha|\leq m}a_{\alpha}(z,\lambda;0)\zeta^{\alpha}.$$

All the results given above are similar in the semiclassical setting.



\section{The Characteristic Set and its Dynamics} 
We want to understand where the operator $P_{\lambda}$ defined by \eqref{psig} is singular and how the singularities of solutions to $P_{\lambda}u=0$ propagate. We see from \eqref{princesymb} that when $\mu\geq0,$  $p$ vanishes only when $\mu=0,$ $\eta=0$.  We follow the notation of \cite{V2}. Let
$$S=\{\mu=0\},\quad N^{*}S\setminus o=\Lambda_{+}\cup\Lambda_{-},\quad \Lambda_{\pm}=N^{*}S\cap\{\pm\xi>0\};$$
so $S\subset X_{-\delta_{0}}$ can be identified with $\partial X.$ We have 
$$N^{*}S=\{(\mu,y,\xi,\eta):\mu=0,\eta=0\}.$$
It follows from \eqref{princesymb} that 
\begin{equation}\label{hamil1}
H_{p}=8\mu\xi\partial_{\mu}+\tilde{H}_{|\eta|_{\mu,y}}-\left(4\xi^2+\partial_{\mu}|\eta|_{\mu,y}^2\right)\partial_{\xi},
\end{equation}
where $\displaystyle \tilde{H}_{|\eta|_{\mu,y}^2}$ is the Hamilton vector field in $T^{*}\partial X,$ and hence 
$$H_{p}|_{N^{*}S}=-4\xi^2\partial_{\xi},$$
this is equivalent to $dp=4\xi^2 d\mu$ at $N^{*}S$, so the characteristic set $\Sigma=\{p=0\}$ is smooth at $N^{*}S$ and $H_{p}$ is radial at $\Lambda_{\pm}.$ 

We work with a compactification of the tangent bundles $T^{*}X_{-\delta_0},$ say $\overline{T^{*}X}_{-\delta_0}.$ The cosphere bundle is $S^{*}X_{-\delta_{0}}=\partial\overline{T^{*}X}_{-\delta_0}$. Let $L_{\pm}=\partial \Lambda_{\pm}$, that is, $\nu\in L_{\pm}$ if $\displaystyle\nu=\lim_{\xi\to\pm\infty}\tilde{\nu},$ where $\tilde{\nu}\in\Lambda_{\pm}$ has $\mu$ and $\eta$ coordinates equal to zero and $\sgn(\xi)=\pm1$. Let $\Sigma_{\pm}=(\Sigma\cap\{\pm\xi>0\})$ be the components of the (classical) characteristic set containing $L_{\pm}$.  Note that it follows from \eqref{princesymb}  that if $\mu>0,$ $\sigma_{2}(P_{\lambda})>0$, so $\Sigma_{\pm}\subset\{\mu\leq0\}.$

In order to understand how the singularities propagate we need to understand the behavior of the flow of the Hamiltonian $H_p$ at the singular points. We have that the equations of the flow of $H_{p}$ are 
\begin{gather*}
\dot{\xi}=-\frac{\partial p}{\partial \mu}=H_{p}\xi=-4\xi^{2}-\partial_{\mu}|\eta|_{\mu,y}^{2}, \quad 
\dot{\mu}=\frac{\partial p}{\partial \xi}=H_{p}\mu=8\mu\xi, \\
\dot{y_{i}}=\frac{\partial p}{\partial \eta_{i}}=H_{p}y_{i}=\partial_{\eta_i}|\eta|_{\mu,y}^2=\sum_{k=1}^{n}(h^{ik}+h^{ki})\eta_{k},\\
\dot{\eta_{i}}=-\frac{\partial p}{\partial y_{i}}=H_{p}\eta_{i}=-\partial_{y_i}|\eta|_{\mu,y}^2=-\sum_{1\leq j,k\leq n}\eta_{j}\eta_{k}\frac{\partial h^{jk}}{\partial y_{i}}.
\end{gather*}

Note that since $h$ is even modulo $O(x^{2k+1}),$ $k\geq 2,$  the coefficients of $H_p$ are Lipschitz and the equations above  have a unique solution for the prescribed initial data. We see that if $\varphi(t)=(\mu(t),y(t),\xi(t),\eta(t))$ solves the equations given above, that is, $\varphi$ is the flow of $H_p$ (also called bicharacteristic), then
\begin{equation}\label{musiglamb}
\mp\dot{\mu}|_{\Sigma_{\pm}}\geq 0,\quad \dot{\xi}|_{\Lambda_{\pm}}\leq 0
\end{equation}

We now analyse the flow at $\Lambda_{\pm}=\{(\mu,y,\xi,\eta):\mu=0,\eta=0,\pm\xi>0\}$. We have that,
\begin{equation}\label{hamieta}
H_{p}|\eta|_{\mu,y}^{2}=8\mu\xi\partial_{\mu}|\eta|_{\mu,y}^2,\quad H_{p}\mu=8\xi\mu,\quad H_{p}|\xi|=-4\text{ }\text{sgn}(\xi)\xi^2+b,\quad  b=-\sgn(\xi)\partial_{\mu}|\eta|_{\mu,y}^2.
\end{equation}
Note that $b$ vanishes at $\Lambda_{\pm}.$

For convenience, to have a better description of the sets $L_{\pm}$ in terms of local coordinates, since they are defined as $\xi\rightarrow\pm\infty$, we rehomogenize in terms of $\displaystyle\hat{\eta}=\eta/|\xi|$ and $\tilde{\rho}=|\xi|^{-1}.$ Then near $N^{*}S$, we can take $\displaystyle \left(|\xi|^{-1},\hat{\eta}\right)$ as projective coordinates on the fibers of the cotangent bundle, with $\tilde{\rho}=|\xi|^{-1}$ defining $S^{*}X_{-\delta_0}$ in $\overline{T^{*}X}_{-\delta_{0}}$. With respect to these new projective coordinate system
$$\partial N^{*}S=\{(\mu,y,\tilde{\rho},\hat{\eta}):\mu=0,\hat{\eta}=0,\tilde{\rho}=0\}.$$

Then, since the metric is assumed to be even modulo $O(x^{2k+1})$, $\displaystyle W=|\xi|^{-1}H_{p}$ is a $C^{k-1}_{\mu},k\geq2,$ vector field and it follows from \eqref{hamieta} that
\begin{equation}
\begin{gathered}
|\xi|^{-1}H_{p}|\hat{\eta}|_{\mu,y}^2=8\mu\frac{\text{sgn}(\xi)}{\xi^2}\partial_{\mu}|\eta|_{\mu,y}^2+2\left[4\xi^2+\partial_{\mu}|\eta|_{\mu,y}^2\right]\text{sgn}(\xi)\frac{|\eta|_{\mu,y}^2}{\xi^4}\\
\Rightarrow |\xi|^{-1}H_{p}|\hat{\eta}|_{\mu,y}^2=8\text{sgn}(\xi)|\hat{\eta}|^2+\tilde{a},
\end{gathered}
\end{equation}
where 
\begin{equation}\label{atil}
\tilde{a}=8\mu\frac{\text{sgn}(\xi)}{\xi^2}\partial_{\mu}|\eta|_{\mu,y}^2+2\partial_{\mu}|\hat{\eta}|_{\mu,y}^2\frac{\sgn(\xi)}{\xi^2}|\hat{\eta}|_{\mu,y}^2.
\end{equation}
Note that $\tilde{a}$ is $C^{k-1}$ and vanishes at $N^{*}S$. Also, it follows from \eqref{hamieta}
\begin{equation}\label{atill}
|\xi|^{-1}H_{p}\mu=8\text{sgn}(\xi)\mu,\quad H_{p}\tilde{\rho}=4\text{sgn}(\xi)+\tilde{a}',\quad \tilde{\rho}=|\xi|^{-1}, \text{ where  }
\tilde{a}'=\frac{\sgn(\xi)}{\xi^2}\partial_{\mu} |\eta|_{\mu,y}^2.
\end{equation}
Note that $\tilde{a}'$ is $C^{k-1}$  and vanishes at $N^{*}S$.

We can see from the equations \eqref{musiglamb}, that if $\phi_{-}(t)=(\mu(t),y(t),\xi(t),\eta(t))$ is a bicharacteristic in $\Sigma_{-}$ ($\sgn(\xi)=-1$) with $\mu(t_0)=-\epsilon_{0}<0$, for some $t_0$, then $\dot{\xi}<0$ and $\dot{\mu}>0,$ so $\xi(t)<0$ and $\mu(t)\rightarrow 0$ as $t\rightarrow\infty$. This then implies that $\phi_{-}(t)$ travels to $L_{-}$ as $t\rightarrow \infty$. Similarly, if $\phi_{+}(t)=(\mu(t),y(t),\xi(t),\eta(t))$ is a bicharacteristics in $\Sigma_{+}$ ($\sgn(\xi)=1$) with $\mu(t_0)=0,$ then since $\mu(t)$ is decreasing in $\Sigma_{+}$, we have that $\mu(t)<0$ as $t\rightarrow \infty$, that is $\phi_{+}(t)$ travels from $L_{+}$ to $\mu=-\epsilon_0$.  

We can summarize the observation given above as: bicharacteristics in $\Sigma_{-}$ $(\sgn\xi=-1)$ travel from $\mu=-\epsilon_0$ to $L_{-},$ while in $\Sigma_{+}$ $(\sgn\xi=1)$ they travel from $L_{+}$ to $\mu=-\epsilon_0.$ Therefore we conclude that $L_{-}=\partial \Lambda_{-}$ is a sink, while $L_{+}$ is a source, in the sense that all nearby bicharacteristics converge to $L_{\pm}$ as the parameter along the bicharacteristic goes to $\mp\infty$.

\begin{figure}[h!]
\begin{center}
\includegraphics[scale=0.70]{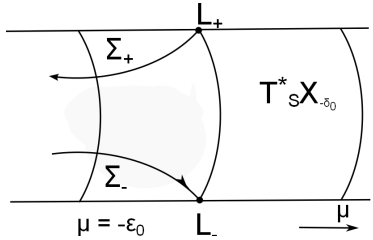}
\caption{The compactified cotangent bundle of $X_{-\delta_0}$ near $S=\{\mu=0\}$, $\overline{T^{*}X}_{-\delta_0}$.}
\end{center}
\end{figure}

We see that $\displaystyle\hat{\rho}_{0}=\widehat{\tilde{p}}+\hat{p}^2, \text{ where } \hat{p}=|\xi|^{-2}p,\text{ }\widehat{\tilde{p}}=|\hat{\eta}|^2,$ is a defining function of $L_{\pm}$ and it satisfies
\begin{equation}\label{inequalityclass}
\text{sgn}(\xi)|\xi|^{-1}H_{p}\hat{\rho}_{0}\geq 16\hat{\rho}_{0}+\mathcal{O}(\hat{\rho}_{0}^{3/2}).
\end{equation}

We will also need to know the form of the principal symbol of $\displaystyle (P_{\lambda}-P_{\lambda}^{*})$. It follows from \eqref{difop1} that the principal symbol of $\displaystyle \frac{1}{2i}(P_{\lambda}-P_{\lambda}^{*})$ at $L_{\pm}$, is given by
\begin{equation}\label{diffprincesymb}
\sigma_{1}\left(\frac{1}{2i}(P_{\lambda}-P_{\lambda}^{*})\right)|_{N^{*}S}=-(4\sgn \xi)|\xi|\text{Im }\lambda=-4(\sgn \xi)\tilde{\rho}^{-1}\im \lambda.
\end{equation}

\section{Semiclassical Dynamics}

We now study the dynamics of the associated semiclassical problem. This is important so that we can understand the high energy behavior, as $|\lambda|\rightarrow \infty.$ We have
$$P_{h,z}=h^{2}P_{h^{-1}z},$$
with $h=|\lambda|^{-1}$ and $z=\lambda/|\lambda|.$ 

We want to work in $\text{Im }\lambda\geq-C$, that is $\text{Im }z\geq -Ch.$ The semiclassical principal symbol of $P_{h,z}$ is given by
\begin{equation}\label{semisymb}
p_{\hbar,z}=\sigma_{2,\hbar}(P_{h,z})=4\mu\xi^2-4(1+a_2)z\xi-z^2+|\eta|_{\mu,y}^2.
\end{equation}

The semiclassical Hamilton vector field is
\begin{equation}\label{hamilsemi}
H_{p_{\hbar,z}}=4(2\mu\xi-(1+a_2)z)\partial_\mu+\widetilde{H}_{|\eta|_{\mu,y}^2}-(4\xi^2-4z\xi\partial_{\mu}a_2+\partial_{\mu}|\eta|_{\mu,y}^2)\partial_{\xi}.
\end{equation}

The real part of \eqref{semisymb} is
\begin{equation}\label{realsymb}
\text{Re }p_{\hbar,z}=4\mu\xi^2-4(1+a_2)(\text{Re }z)\xi-\left[(\text{Re }z)^2-(\text{Im }z)^2\right]+|\eta|_{\mu,y}^2,
\end{equation}
and the imaginary part
\begin{equation}\label{imasemi}
\text{Im }p_{\hbar,z}=-2(\text{Im }z)\left[2(1+a_2)\xi+\text{Re }z\right].
\end{equation}

Therefore, if \eqref{imasemi} vanishes, for $z$ non-real ($\text{Im }z\neq0$), we have
\begin{equation}\label{semire}
\text{Re }p_{\hbar,z}=(\text{Re }z)^2\left[(1+a_2)^{-2}\mu+1\right]+(\text{Im }z)^2+|\eta|_{\mu,y}^2>0, \quad \text{for }\mu>0.
\end{equation}
Then $p_{\hbar,z}$ is semiclassically elliptic on $T^{*}X_{-\delta_0}$, but not at $S^{*}X_{-\delta_0}$, we will explain this bellow. 

As in the classical setting, the characteristic set can be divided into two components $\Sigma_{\hbar,\pm}$, where $\Sigma_{\hbar}=\{p_{\hbar,z}=0\}$. As we saw above, if $z$ is non-real, $p_{\hbar,z}=0$ implies
$$2(1+a_2)\xi+\text{Re }z=0,$$
so it follows from \eqref{semire} that $\text{Re }p_{\hbar,z}$ does not vanish, hence
$$\Sigma_{\hbar}=\Sigma_{\hbar,+}\cup\Sigma_{\hbar,-},\quad \Sigma_{\hbar,\pm}=\Sigma_{\hbar}\cap\{\pm(2(1+a_2)\xi+\text{Re }z)>0\}.$$

If we introduce $\displaystyle(\mu,y,\rho,\hat{\eta}),\text{ }\rho=|\xi|^{-1},\text{ }\hat{\eta}=\eta/|\xi|,$ as projective coordinates in a neighborhood of $L_{\pm}$ in $\overline{T}^{*}X_{-\delta_0},$ we have
$$\rho^2p_{\hbar,z}=4\mu-4(1+a_2)(\sgn\xi)z\rho-z^2\rho^2+|\hat{\eta}|_{\mu,y}^2.$$
Hence $\displaystyle\rho^2\text{Im }p_{\hbar,z}=-2\text{Im }z\left[2(1+a_2)(\sgn\xi)\rho+\rho^2\text{Re }z\right],$ which vanishes when $\rho=0$, that is, at $S^{*}X_{-\delta_0}.$ 

Then at $\rho=0$ when $\text{Im }z$ is bounded away from $0$, we analyse the flow of $H_{p_{\hbar,z}}$ similarly as we did in the classical setting and obtaing:
\begin{equation}
\begin{gathered}
\rho\text{ small },\text{ Im }z\geq 0\Rightarrow(\sgn\xi)\text{Im }p_{\hbar,z}\leq0\Rightarrow \pm\text{Im }p_{\hbar,z}\leq 0\text{ near }\Sigma_{\hbar,\pm}\\
\rho\text{ small },\text{ Im }z\leq 0\Rightarrow(\sgn\xi)\text{Im }p_{\hbar,z}\geq0\Rightarrow \pm\text{Im }p_{\hbar,z}\geq 0\text{ near }\Sigma_{\hbar,\pm}
\end{gathered}
\end{equation}

This implies, as in the classical case, that for $\text{Im }z>0$, $L_{+}$ is a source and $L_{-}$ is a sink, while for $P_{h,z}^{*}$ the directions are reversed. The directions are also reversed if $\text{Im }z$ switches sign. 

The following lemma was proved in \cite{V2}.

\begin{lemma}\label{lemmavasy}
There exists $\epsilon_0>0$ such that the following holds. All semiclassical null-bicharacteristics in $(\Sigma_{\hbar,+}\setminus L_{+})\cap\{-\epsilon_0\leq \mu\leq \epsilon_0\}$ go to either $L_{+}$ or to $\mu=\epsilon_0$ in the backward direction and to $\mu=\epsilon_0$ or $\mu=-\epsilon_0$ in the forward direction, while all semiclassical null-bicharacteristics in $(\Sigma_{\hbar,-}\setminus L_{-})\cap \{-\epsilon_0\leq\mu\leq\epsilon_0\}$ go to $L_{-}$ or to $\mu=\epsilon_0$ in the forward direction and to $\mu=\epsilon_0$ or $\mu=-\epsilon_0$ in the backward direction. 

For $\text{Re }z>0$, only $\Sigma_{\hbar,+}$ enters $\mu>0$, so the $\mu=\epsilon_0$ possibility only applies to $\Sigma_{\hbar,+}$, while for $\text{Re }z<0$, the analogous remark to $\Sigma_{\hbar,-}.$
\end{lemma}

We will consider the lemma \ref{lemmavasy} above when $g$ is a non-trapping metric, that is, every geodesic $\gamma(t)$ reaches $\partial X$ as $ t\rightarrow \pm\infty$, then the case $\mu=\epsilon_0$ does not need to be considered. We have a stronger result: for sufficiently small $\epsilon_0$, and for $\text{Re }z>0$, any bicharacteristics in $\Sigma_{\hbar,+}$ in $-\epsilon_0\leq \mu$ goes to $L_{+}$ at $-\infty$ (backward direction), and to $\mu=-\epsilon_0$ at $\infty$ (forward direction), except the constant bicharacteristics at $L_{+}$; as for $\Sigma_{\hbar,-}$, all bicharacteristics in $-\epsilon_0\leq \mu$ stay in $-\epsilon_0\leq\mu\leq0$, and go to $L_{-}$ in the forward direction and to $\mu=-\epsilon_0$ in the backward direction, except the constant bicharacteristics at $L_{-}$. 






\section{The Complex Absorption Operator}
Finally we now want to work in a compact manifold without boundary. We consider a compact manifold without boundary $X_{\mu_0}=\{\mu>\mu_{0}\},$ $\mu_0=-\epsilon_0<0,$ with $\epsilon_{0}>0$ as above.

To deal with the region where the operator $P_{\lambda}$ is not elliptic, we introduce a complex absorbing operator $Q_{\lambda}\in\Psi^{2}(X)$ with real principal symbol $q$ supported in, say, $\mu<-\epsilon_1,$ where its Schwartz kernel is also supported in the corresponding region,  such that $p\pm iq$ is elliptic in a neighborhood of $\partial X_{\mu_0}$, and that $\pm q\geq 0$ holds on $\Sigma_{\pm}$. We can see that such a $q$ exists since the characteristic set of $p$ has two disjoint connected components, $\Sigma_{\pm}$, and $p$ is elliptic away from $\Sigma_{\pm}$, therefore $p\pm iq$ is also elliptic, independently of the choice of $q$; we just need to make the support of $q$ in a small neighborhood of $\Sigma_{\pm}$, $q$ elliptic on $\Sigma_{\pm}$ at $\mu=-\epsilon_{0}$, with the same sign in $\Sigma_{\pm}$ as near $\{\mu=0\}$. We then extend $P_{\lambda}$ and $Q_{\lambda}$ to $X$ so that $p\pm iq$ are elliptic near $\partial X_{\mu_0}$. We then take $Q_{\lambda}$ to be any operator whose principal symbol is $q$. 

Other conditions that are also satisfied by $Q_{\lambda}$ are the following: for $\nu\in S^{*}X\cap\Sigma$, let $\gamma_{+}(\nu)$, respectively $\gamma_{-}(\nu)$ denote the image of the forward, resp. backward, half-bicharacteristics of $p$ from $\nu$. We say that $\gamma_{\pm}(\nu)$ tends to $L_{\pm}$ ($\gamma_{\pm}(\nu)\rightarrow L_{\pm}$) if for any neighborhood $O$ of $L_{\pm}$, $\gamma_{\pm}(\nu)\cap O\neq\emptyset$. From what we have seen in the previous section, if we denote by $\text{ell}(Q_{\lambda})$ the elliptic set of $Q_{\lambda}$, it follows from the assumptions on the operator $Q_{\lambda}$, namely its principal symbol has support in $\mu<-\epsilon_1$, for some $\epsilon_1>0,$ $\pm q\geq 0$ on $\Sigma_{\pm}$ and $\Sigma_{\pm}\subset \text{ell}(Q_{\lambda}),$ we have that
\begin{equation}\label{ellq}
\begin{gathered}
\nu\in \Sigma_{-}\setminus L_{-}\Rightarrow \gamma_{+}(\nu)\rightarrow L_{-} \text{ and }\gamma_{-}(\nu)\cap\text{ell}(Q_{\lambda})\neq\emptyset,\\
\nu\in \Sigma_{+}\setminus L_{+}\Rightarrow \gamma_{-}(\nu)\rightarrow L_{+} \text{ and }\gamma_{+}(\nu)\cap\text{ell}(Q_{\lambda})\neq\emptyset.
\end{gathered}
\end{equation}

Indeed we have shown previously that if $\displaystyle \nu\in \Sigma_{\pm}\setminus L_{\pm}$, then $\displaystyle\gamma_{\mp}(\nu)\rightarrow L_{\pm}$, but since $L_{\pm}\subset\Sigma_{\pm}\subset \text{ell}(Q_{\lambda})$, $\gamma_{\mp}(\nu)\cap\text{ell}(Q_{\lambda})$ is not empty. That is, all forward and backward half-(null)bicharacteristics of $P_{\lambda}$ either enter the elliptic set of $Q_{\lambda}$, or go to $\Lambda_{\pm},$ i.e., $L_{\pm}$ in $S^{*}X.$ 

For the semiclassical problem, when $z$ is non-real, $Q_{\lambda}$ needs to have more properties, and that depends if $g$ is non-trapping. We then choose $Q_{\hbar,z}=h^{2}Q_{h^{-1}z}\in\Psi^{2}_{\hbar}(X)$ with semiclassical principal symbol $q_{\hbar,z}$ to have all the above requirements and be semicassically elliptic near $\mu=\mu_0,$ that is, $p_{\hbar,z}-iq_{\hbar,z}$ and its complex conjugate are elliptic near $\{\mu=\mu_0\},$ and $\pm q_{\hbar,z}\geq 0$ holds for $z$ real on $\Sigma_{\hbar,\pm}$. In \cite{V2} Vasy constructed an example of an operator $Q_{\hbar,z}$ which satisfies the above hypothesis and it is a holomorphic family at least in a strip.


By \eqref{imasemi} we can take (with $C>0$)
\begin{equation}
q_{\hbar,z}=2\left[2(1+a_2)\xi+z\right](\xi^2+|\eta|_{\mu,y}^2+z^2+C^2h^2)^{1/2}\chi(\mu),
\end{equation}
where $\chi\geq 0$ has support contained in a neighborhood of $\mu_0.$ 

The corresponding semiclassical principal symbol is
$$\sigma_{2,\hbar}(Q_{\hbar,z})=2\left[2(1+a_2)\xi+z\right](\xi^2+|\eta|_{\mu,y}^2+z^2)^{1/2}\chi(\mu),$$
and we take $Q_{\hbar,z}$ as the quantization of this symbol. 

Here we assume that the square root of a complex number is defined in $\mathbb{C}\setminus[0,-\infty)$, so $q_{\hbar,z}$ is defined away from $h^{-1}z\in \pm i[C,+\infty). $ We then see that the symbol $\xi^2+|\eta|_{\mu,y}^2+\lambda^2$ is elliptic in $(\xi,\eta,\text{Re }\lambda,\text{Im }\lambda)$ if $|\text{Im }\lambda|<C'|\text{Re }\lambda|$, for some $C'>0.$ 

Therefore, assuming $g$ is even modulus $O(x^{2k+1}),k\geq 2$, it follows from Proposition \ref{ellipticregu} $\displaystyle p(\mu,y,\xi,\eta)\in S_{st}^{2}(k),$ and if $\displaystyle -k<\delta<s\leq k$, we have that
$$P_{\lambda}-iQ_{\lambda}\text{ elliptic at }\nu,\text{ } \nu\not\in [\text{WF}^{s-2}((P_{\lambda}-iQ_{\lambda})u)\cup\text{WF}^{\delta}(u)]\Rightarrow \nu\not\in\text{WF}^{s}(u).$$
In particular, if $\displaystyle(P_{\lambda}-iQ_{\lambda})u\in H^{s-2}$, $u\in H^{\delta}$ and $p-iq$ is elliptic at $\nu$, then $\nu\not\in\text{WF}^{s}(u).$ Analogous conclusions apply to $P^{*}_{\lambda}+iQ_{\lambda}^{*};$ since both $p$ and $q$ are real, $p-iq$ is elliptic if and only if $p+iq$ is. 


Moreover, if we let $\tilde{\gamma}_{\pm}(\nu)$ be the forward (+) or backward (-) bicharacteristic from $\nu$, then it follows from \eqref{ellq} and \eqref{ellipticregu} that, for $-k<s<k+1$,
\begin{equation}
\begin{gathered}
\nu\not\in WF^{s}(u)\text{ and }WF^{s-1}((P_{\lambda}-iQ_{\lambda})u)\cap\tilde{\gamma}_{\pm}(\nu)=\emptyset\Rightarrow \tilde{\gamma}_{\pm}(\nu)\cap WF^{s}(u)=\emptyset,\\
\nu\not\in WF^{s}(u)\text{ and }WF^{s-1}((P_{\lambda}^{*}+iQ_{\lambda}^{*})u)\cap\tilde{\gamma}_{\pm}(\nu)=\emptyset\Rightarrow \tilde{\gamma}_{\pm}(\nu)\cap WF^{s}(u)=\emptyset.
\end{gathered}
\end{equation} 

We actually have that in the region where $q\neq0$, by elliptic regularity $\nu\not\in WF^{s-1}((P_{\lambda}-iQ_{\lambda})u)$ implies $\nu\not\in WF^{s+1}(u)$. We then treat the case when $q=0$ in the following section.

\section{Meromorphic Extension of the Resolvent}

We will show below, Theorems \ref{meroextention} and \ref{thmpq}, that all the desired properties that we want to show for $R(\lambda)=\left(\Delta_{g}-\frac{n^2}{4}-\lambda^2\right)^{-1}$ are satisfied by the inverse of the operator $P_{\lambda}-iQ_{\lambda}$. Now we will show how to transfer those properties from $(P_{\lambda}-iQ_{\lambda})^{-1}$ to $R(\lambda).$ Here we will assume that the operator $P_{\lambda}-iQ_{\lambda}$ acts on $H^{s}(X_{-\delta_0})$, for some $s\in\mathbb{R}$. Note that $X_{-\delta_0}$ is a pre-compact manifold without boundary and $X\subset X_{-\delta_0}$ can be seen as a compact subset of $X_{-\delta_0}$.  

First, since the Laplace-Beltrami operator $\Delta_g$ is self-adjoint, for any $f\in \mathcal{\dot{S}}^{\infty}(X)$, the problem
\begin{equation}\label{wavevasy}
\left(\Delta_g-\frac{n^2}{4}-\lambda^2\right)u=f
\end{equation}
has a unique solution $u=R(\lambda)f\in L^2(X,|dg|),$ for $\text{Im }\lambda>>0.$ 

Now let $\varphi$ be such that $\displaystyle e^{\phi}=\mu^{1/2}(1+\mu)^{-1/4}$ near $\mu=0$ with $\displaystyle 4\mu^2(\partial_{\mu}\varphi)^2<1$. Let $\tilde{f}_{0}\in\mathcal{\dot{S}}^{\infty}(X)$ be such that $\displaystyle \tilde{f}_{0}=e^{i\lambda\varphi}x^{-(n+2)/2}f$ in $\mu\geq0$. Let $\tilde{f}$ be any smooth extension of $\tilde{f}_{0}$ to $X_{-\delta_0}=((-\delta_0,0)_{\mu}\times \partial X)\cup X,$ where $\delta_0>0$ is such that $P_{\lambda}-iQ_{\lambda}$ is defined in $X_{-\delta_0}$. 

Then by Theorem \ref{thmpq} below, we can take $\tilde{u}_{0}=(P_{\lambda}-iQ_{\lambda})^{-1}\tilde{f}$. We have that, since $\tilde{f}\in C^{\infty}(X)$, $\tilde{u}_{0}\in C^{\infty}(X)$ (actually $f\in H^{s}$ for some $s$, so $\tilde{u}_{0}\in H^{s}$, but iterating the statement of Theorem \ref{meroextention}, we can deduce that $\tilde{u}_{0}\in C^{\infty}(X)).$ We have that $\displaystyle(P_{\lambda}-iQ_{\lambda})\tilde{u}_{0}=\tilde{f}.$

Now since $Q_{\lambda}$ is supported in $\mu\leq0$, if we let $\tilde{u}=e^{-i\lambda\varphi}x^{(n+2)/2}x^{-1}\tilde{u}_{0}|_{\mu>0}$, then $\tilde{u}\in x^{n/2}e^{-i\lambda\varphi}C^{\infty}(X)$, and
$$P_{\lambda}e^{i\lambda\varphi}x^{-(n+2)/2}x\tilde{u}=\tilde{f}|_{\mu>0}=\tilde{f}_{0}=e^{i\lambda\varphi}x^{-(n+2)/2}f,$$
but since 
$$P_{\lambda}=e^{i\lambda\varphi}\mu^{-(n+2)/4-1/2}\left(\Delta_g-\frac{n^2}{4}-\lambda^2\right)e^{-i\lambda\varphi}\mu^{(n+2)/4-1/2},$$
and $x=\mu^{1/2},$ we have 
$$\left(\Delta_g-\frac{n^2}{4}-\lambda^2\right)\tilde{u}=f.$$

Since \eqref{wavevasy} has a unique solution and $\tilde{u}\in L^{2}(X,|dg|)$ for $\mu>0$, $u=\tilde{u}.$ We then have that for $f\in \mathcal{\dot{S}}^{\infty}(X),$

\begin{equation}\label{resolpq}
R(\lambda)f=e^{-i\lambda\varphi}x^{(n+2)/2-1}(P_{\lambda}-iQ_{\lambda})^{-1}e^{i\lambda\varphi}x^{-(n+2)/2-1}f.
\end{equation} 

We have proved above that the formula \eqref{resolpq} is independent of the choice of $Q_{\lambda}$. Then Theorem \ref{mainthm} follows directly from \eqref{resolpq} and Theorems \ref{meroextention}, \ref{thmpq}, proved bellow. 

We are going to show the meromorphy of the inverse of the operator $P_{\lambda}-iQ_{\lambda}$ by applying analytic Fredholm theory, which will follow from regularity results, namely propagation of singularities, and the properties of the complex absorping operator. We now analyze the propagation of the singularities near $\Lambda_{\pm}.$ We will prove the following theorem.

\begin{thm}\label{meroextention}
Let $P_{\lambda},$ $Q_{\lambda}$ be as above. Let $-k<m<s\leq k+1$. Then 
$$\text{if } U_{m,s}=\{u\in H^{m}(X_{-\delta_0}):(P_{\lambda}-iQ_{\lambda})u\in H^{s-1}(X_{-\delta_0})\}\cap H^{s}(X_{-\delta_0}),$$
$$P_{\lambda}-iQ_{\lambda}:U_{m,s}\rightarrow H^{s-1}(X_{-\delta_0})$$
is an analytic family of Fredholm operators on
\begin{equation}
\mathbb{C}_{s}=\{\lambda\in\mathbb{C}:\text{ Im }\lambda>1/2-s\}.
\end{equation}
\end{thm} 

We start by proving our main regularity result. We prove an analogue of Vasy's propagation of singularities result proved in \cite{V2}. The proof is adapted to the case of pseudodifferential operators with rough symbols. Since the regularity of the operator $P_{\lambda}$ depends on the order of evenness of the metric $g$, we have the following.

\begin{prop}\label{propsing} Suppose $s> m>1/2-\text{Im }\lambda$, $\displaystyle m< s<k+1$, and $\text{WF}^{m}(u)\cap\Lambda_{\pm}=\emptyset.$ Then
$$\Lambda_{\pm}\cap\text{WF}^{s-1}(P_{\lambda}u)=\emptyset\Rightarrow \Lambda_{\pm}\cap\text{WF}^{s}(u)=\emptyset.$$
\end{prop} 
\textbf{Proof.} The proof is based on a positive commutator estimate. We will consider $C_{\epsilon}^{*}C_{\epsilon}$ with $C_{\epsilon}\in\Psi^{s-1/2-\delta}(X)$ for $\epsilon>0$, uniformly bounded in $\Psi^{s-1/2}(X)$ as $\epsilon\rightarrow 0,$ where $\delta=s-m$. We let
$$c=\phi(\rho_{0})\tilde{\rho}^{-s+1/2},\quad c_{\epsilon}=c(1+\epsilon\tilde{\rho}^{-1})^{-\delta},$$
where
$$\tilde{\rho}=|\xi|^{-1},\quad \rho_{0}=\frac{|\eta|^2}{|\xi|^2}+|\xi|^{-4},$$ 
and $\phi\in C_{c}^{\infty}(\mathbb{R})$ is identically 1 near 0, $\phi'\leq 0$ on $(0,\infty)$ and $\phi$ is supported sufficiently close to 0 so that  
\begin{equation}\label{inephi}
\text{if }(\xi,\eta)\text{ is such that }\rho_{0}(\xi,\eta)\in\text{supp }\phi'\Rightarrow \sgn(\xi)|\xi|^{-1}H_{p}\rho_{0}>0;
\end{equation}
such $\phi$ exists since 
$$(\text{sgn }\xi)|\xi|^{-1}H_{p}\rho_{0}\geq 8\rho_{0}+O(\rho_{0}^{3/2}).$$
We choose $\phi$ so that $\sqrt{-\phi\phi'}$ is $C^{\infty}$ on $[0,\infty)$. Note that the sign of $H_{p}\tilde{\rho}^{-s+1/2}$ depends on the sign of $-s+1/2$.

Now let $C\in\Psi^{s-1/2}(X)$ have principal symbol $c$, and have $\text{WF}'(C)\subset \text{supp }\phi\circ\rho_{0}$, and let $C_{\epsilon}=CS_{\epsilon},$ where $S_{\epsilon}\in\Psi^{-\delta}(X)$ has principal symbol $\displaystyle (1+\epsilon\tilde{\rho}^{-1})^{-\delta}=(1+\epsilon|\xi|)^{-\delta},$ which is uniformly bounded in $\Psi^{0}(X)$ for $\epsilon>0$, converging to Id in $\Psi^{\delta'}(X)$ for $\delta'>0$ as $\epsilon\rightarrow 0.$ Then the principal symbol of $C_{\epsilon}$ is $c_{\epsilon}\in S^{s-1/2-\delta}$. We have
\begin{equation}\label{positivecomu}
P_{\lambda}^{*}C_{\epsilon}^{*}C_{\epsilon}-C_{\epsilon}^{*}C_{\epsilon}P_{\lambda}=(P_{\lambda}^{*}-P_{\lambda})C_{\epsilon}^{*}C_{\epsilon}+[P_{\lambda},C_{\epsilon}^{*}C_{\epsilon}].
\end{equation}

We see that by Lemmas \ref{compoab} and \ref{comporem}, since the symbol of the operator $P_{\lambda}$ is in $S^{2}_{st}(k),$ the operator $\displaystyle C_{\epsilon}^{*}C_{\epsilon}P_{\lambda}$ is well defined, bounded in $H^{r+2s+1-2\delta}$ if $\displaystyle 2s-1-2\delta-k<r\leq k+1+2\delta-2s,$ 
and its full symbol $c_{\epsilon}^2\circ p$ satisfies
\begin{equation}\label{resto}
r(\mu,y,\xi,\eta)=c_{\epsilon}^2\circ p(\mu,y,\xi,\eta)-\sum_{|\alpha|< 2s-1-2\delta}(i^{-|\alpha|}/\alpha!)\partial_{\xi}^{\alpha_1}\partial_{\eta}^{\alpha'}c_{\epsilon}^2(\mu,y,\xi,\eta)\partial_{y}^{\alpha'}\partial_{\mu}^{\alpha_1}p(\mu,y,\xi,\eta),
\end{equation}
with $r(\mu,y,\xi,\eta)\in S^{2}\left(k+1-2s+2\delta\right)$. We can then compute the principal symbol of (\ref{positivecomu}), which is given by
$$\sigma_{2s}\left(i(P_{\lambda}^{*}C_{\epsilon}^{*}C_{\epsilon}-C_{\epsilon}^{*}C_{\epsilon}P_{\lambda})\right)=\sigma_{1}(i(P_{\lambda}^{*}-P_{\lambda}))c_{\epsilon}^2+H_{p}(c_{\epsilon}^2)=\sigma_{1}(i(P_{\lambda}^{*}-P_{\lambda}))c_{\epsilon}^2+2c_{\epsilon}H_{p}c_{\epsilon},$$
since the principal symbol of $C_{\epsilon}$ and $C_{\epsilon}^{*}$ are the same, $c_{\epsilon}$, and the principal symbol of $[P_{\lambda},C_{\epsilon}^{*}C_{\epsilon}]$ is $-iH_{p}c_{\epsilon}^2.$\\
It follows from \eqref{diffprincesymb} that
\begin{equation}
\sigma_{1}\left(\frac{1}{2i}(P_{\lambda}-P_{\lambda}^{*})\right)|_{\mu=0}=-4(\text{Im }\lambda)\xi.
\end{equation}
We have then 
$$\sigma_{1}(i(P_{\lambda}^{*}-P_{\lambda}))c_{\epsilon}^2+2c_{\epsilon}H_{p}c_{\epsilon}$$
$$= 8\sgn\xi\left(-\text{Im }\lambda\phi+\left(-s+\frac{1}{2}\right)\phi+\frac{\sgn\xi}{4}(\tilde{\rho}H_{p}\rho_{0})\phi'+\delta\frac{\epsilon}{\tilde{\rho}+\epsilon}\phi\right)\phi\tilde{\rho}^{-2s}(1+\epsilon\tilde{\rho}^{-1})^{-2\delta}.$$
Thus,
\begin{equation}
\begin{gathered}
(\sgn\xi)\sigma_{2s+1-2\delta}(i(P_{\lambda}^{*}C_{\epsilon}^{*}C_{\epsilon}-C_{\epsilon}^{*}C_{\epsilon}P_{\lambda}))
\leq \\ -8\left(s-\frac{1}{2}+\text{Im }\lambda-\delta\right)\tilde{\rho}^{-2s}(1+\epsilon\tilde{\rho}^{-1})^{-2\delta}\phi^2
+2((\sgn\xi)\tilde{\rho}H_{p}\rho_{0})\tilde{\rho}^{-2s}(1+\epsilon\tilde{\rho}^{-1})^{-2\delta}\phi'\phi.
\end{gathered}
\end{equation}
We see that the first term on the right hand side is negative if $s-1/2+\text{Im }\sigma-\delta>0$ and this has the same sign of the term with $\phi'$. We then let 
\begin{equation}
\begin{gathered}
B\in\Psi^{s}\text{ with }\sigma_{s}(B)=\sqrt{8}\left(\sqrt{s-\frac{1}{2}+\text{Im }\lambda}\right)\phi(\rho_0)\tilde{\rho}^{-s},\\
B_{2,\epsilon}\in\Psi^{s}\text{ with }\sigma_{s}(B_2,\epsilon)=\sqrt{8}\sqrt{\frac{\delta\epsilon}{\tilde{\rho}+\epsilon}}\phi(\rho_0)\tilde{\rho}^{-s},\\
F_{\epsilon}\in\Psi^{2s-1}(k+1-2s)\text{ with }\sigma_{2s-1}(F_{\epsilon})=2(\sgn \xi)(\tilde{\rho}H_{p}\rho_0)\tilde{\rho}^{-2s}(1+\epsilon\tilde{\rho}^{-1})^{-2\delta}\phi'(\rho_0)\phi(\rho_0).
\end{gathered}
\end{equation}  
Note that $\displaystyle B_{2,\epsilon}$ is uniformly bounded in $\Psi^{s}(X)$ as $\epsilon\rightarrow 0$ and $F_{\epsilon}$ is uniformly bounded in $\Psi^{2s-1}(k+1-2s)$. It also follows from our choice of $\phi$ that $\displaystyle (\sgn \xi)(\tilde{\rho}H_{p}\rho_0)\phi'(\rho_0)\phi(\rho_0)\leq 0$, so $\sigma_{2s-1}(F_{\epsilon})\leq 0$. We then write
$$(\sgn \xi) i(P_{\lambda}^{*}C_{\epsilon}^{*}C_{\epsilon}-C_{\epsilon}^{*}C_{\epsilon}P_{\lambda})=-S_{\epsilon}^{*}(B^{*}B-B_{2,\epsilon}^{*}B_{2,\epsilon})S_{\epsilon}+F_{\epsilon}.$$

We have that

\begin{equation}\label{mel1}
\left\langle i(P_{\lambda}^{*}C_{\epsilon}^{*}C_{\epsilon}-C_{\epsilon}^{*}C_{\epsilon}P_{\lambda})u,u\right\rangle=
2\text{Im}\left\langle C_{\epsilon}^{*}C_{\epsilon}u,P_{\lambda}u\right\rangle\leq 2\left|\left\langle C_{\epsilon}^{*}C_{\epsilon}u,P_{\lambda}u\right\rangle\right|.
\end{equation}
Using Cauchy inequality
$$
\left|\left\langle C_{\epsilon}^{*}C_{\epsilon}u,P_{\lambda}u\right\rangle\right|\leq C\|C_{\epsilon}u\|_{H^{1/2}}^2\|C_{\epsilon}P_{\lambda}u\|_{H^{-1/2}}^2,
$$
but since the symbol of $S_{\epsilon}B$ is an elliptic multiple of $\displaystyle\tilde{\rho}^{-\frac{1}{2}}\sigma(C_{\epsilon})$, if we take $Y\in\Psi^{1/2}$ such that $\sigma(Y)=\tilde{\rho}^{-1/2}$, $\sigma(BS_{\epsilon})=\sigma(YC_{\epsilon})$ it leads to
\begin{equation}\label{mel2}
\left|\left\langle C_{\epsilon}^{*}C_{\epsilon}u,P_{\lambda}u\right\rangle\right|\leq \epsilon^2\|BS_{\epsilon}u\|_{L^2}^2+C\epsilon^{-2}\|BS_{\epsilon}P_{\lambda}u\|_{H^{-1}}^2
\end{equation}
We also have,
$$-\|BS_{\epsilon}u\|_{L^2}^{2}+\|B_{2,\epsilon}S_{\epsilon}u\|_{L^2}^{2}+ \left\langle F_{\epsilon}u,u\right\rangle=2(\sgn \xi)\text{Im}\left\langle C_{\epsilon}^{*}C_{\epsilon}u,P_{\lambda}u\right\rangle.$$
Thus, if we take $\epsilon$ small enough,
$$
\|BS_{\epsilon}u\|_{L^2}^2\leq C\left|\left\langle C_{\epsilon}^{*}C_{\epsilon}u,P_{\lambda}u\right\rangle\right|+\left\langle F_{\epsilon}u,u\right\rangle.
$$
Then using (\ref{mel2}), we have
\begin{equation}\label{mel3}
\|BS_{\epsilon}u\|_{L^2}^2\leq C\|BS_{\epsilon}P_{\lambda}u\|_{H^{-1}}^2+\left\langle F_{\epsilon}u,u\right\rangle.
\end{equation}
Now it follows from the sharp G{\aa}rding inequality (proposition (\ref{Garding})) that 
\begin{equation}\label{mel4}
\left\langle F_{\epsilon}u,u\right\rangle\leq C\|u\|_{H^{s-\delta}}^{2},
\end{equation}
provided $k+1-s+\delta>0$ and $\displaystyle -1\leq \frac{2k+2-2s+2\delta}{k+3-s+\delta}.$

So it follows from (\eqref{mel3}) and (\eqref{mel4})  
\begin{equation}
\|BS_{\epsilon}u\|_{L^2}^{2}\leq C\left(\|BS_{\epsilon}P_{\lambda}u\|_{H^{-1}}^2+\|u\|_{H^{s-\delta}}^2\right).
\end{equation}
This shows that $u$ is in $H^{s}$ in the elliptic set of $B$, provided that $u$ is microlocally in $H^{s-\delta}$ and $P_{\lambda}u$ is in $H^{s-1}$ on the elliptic set of $B$, which contains $\Lambda_{\pm}$. Finally using induction, starting with $s-\delta=m$ and improving regularity in each step by $\leq 1/2$ gives what we want to prove.
$\Box$

Since when taking $P_{\lambda}^{*}$ in place of $P_{\lambda}$ the principle symbol of $\displaystyle \frac{1}{2i}(P_{\lambda}-P_{\lambda}^{*})$ switches sign, we have the following analogous proposition.
\begin{prop}\label{propsinga}
For $s<1/2+\text{Im }\lambda,$ $s<k+1,$ and $O$ a neighborhood of $\Lambda_{\pm}$,
$$\text{WF}^{s}(u)\cap(O\setminus\Lambda_{\pm})=\emptyset,\quad \text{WF}^{s-1}(P_{\lambda}^{*}u)\cap\Lambda_{\pm}=\emptyset\Rightarrow \text{WF}^{s}(u)\cap\Lambda_{\pm}=\emptyset.$$ 
\end{prop}

\textbf{Proof.} We now apply $P_{\lambda}$ in place of $P_{\lambda}^{*}$. Consider the operator $C_{\epsilon}$ as in the proof of Proposition \ref{propsing}. We now consider 
$$P_{\lambda}C_{\epsilon}^{*}C_{\epsilon}-C_{\epsilon}^{*}C_{\epsilon}P_{\lambda}^{*}=(P_{\lambda}-P_{\lambda}^{*})C_{\epsilon}^{*}C_{\epsilon}+[P_{\lambda}^{*},C_{\epsilon}^{*}C_{\epsilon}].$$
We see that the operators above are well defined, the symbol of $C_{\epsilon}^{*}C_{\epsilon}P_{\lambda}^{*}$ is in $S^{2s+1-2\delta}(k+1-2s+2\delta).$ Then
\begin{equation}
\begin{gathered}
\sigma_{2s}(i(P_{\lambda}C_{\epsilon}^{*}C_{\epsilon}-C_{\epsilon}^{*}C_{\epsilon}P_{\lambda}^{*}))=\sigma_{1}(i(P_{\lambda}-P_{\lambda}^{*}))c_{\epsilon}^2+2c_{\epsilon}H_{p}c_{\epsilon}\\
=8\sgn\xi\left((\text{Im }\lambda)\phi+\left(-s+\frac{1}{2}\right)\phi+\frac{\sgn\xi}{4}(\tilde{\rho}H_{p}\rho_{0})\phi'+\delta\frac{\epsilon}{\tilde{\rho}+\epsilon}\phi\right)\phi\tilde{\rho}^{-2s}(1+\epsilon\tilde{\rho}^{-1})^{-2\delta}.
\end{gathered}
\end{equation}
Since we want it to be positive, we need $\displaystyle s-1/2-\text{Im }\lambda<0,$ that is $s<1/2+\text{Im }\lambda$. On the other hand, the term with $\phi'$ now has the opposite sign. We then make the assumption on the support of $\phi'$ so that this is in $O\setminus\Lambda_{\pm},$ that is $\phi'(|\xi|)=0$ if $\xi\not\in \pi_{\xi}(O\setminus\Lambda_{\pm})$, where $\pi_{\xi}$ is the projection onto the $\xi$ coordinate. This is archieved if we make the support of $\phi$ small enough and $\phi$ still indentically $1$ near $0$. Here we take
\begin{equation}
\begin{gathered}
B\in\Psi^{s}\text{ with }\sigma_{s}(B)=\sqrt{8}\sqrt{\left(\frac{1}{2}+\text{Im }\lambda-s\right)}\phi(\rho_0)\tilde{\rho}^{-s},
B_{2,\epsilon}\in\Psi^{s}\text{ with }\sigma_{s}(B_2,\epsilon)=\sqrt{8}\sqrt{\frac{\delta\epsilon}{\tilde{\rho}+\epsilon}}\phi(\rho_0)\tilde{\rho}^{-s}\\
\end{gathered}
\end{equation}  
and $F_{\epsilon}$ as in the proof of proposition \ref{propsing}. Now we have
\begin{equation}
(\sgn \xi)i(P_{\lambda}C_{\epsilon}^{*}C_{\epsilon}-C_{\epsilon}^{*}C_{\epsilon}P_{\lambda}^{*}))=S_{\epsilon}^{*}(B^{*}B+B_{2,\epsilon}^{*}B_{2,\epsilon})S_{\epsilon}+F_{\epsilon}.
\end{equation}
Then
$$
\|BS_{\epsilon}u\|_{L^2}^2\leq \|BS_{\epsilon}u\|_{L^2}^2+\|B_{2,\epsilon}S_{\epsilon}^2\|_{L^2}^2=-\langle F_{\epsilon}u,u\rangle+2(\sgn \xi)\text{Im}\langle C_{\epsilon}^{*}C_{\epsilon}u,P_{\lambda}^{*}u\rangle.
$$
If we follow as in the proof of proposition \ref{propsing}, use the hypothesis that $\displaystyle \text{WF}^{s}(u)\cap\left(O\setminus\Lambda_{\pm}\right)=\emptyset$ and $\Sigma(F_{\epsilon})\subset O\setminus\Lambda_{\pm},$ then using Cauchy-Schwarz inequality we can control $-\langle F_{\epsilon}u,u\rangle$ by $C(\|Au\|_{H^{s}}+\|u\|_{H^{-N}}),$ for any $N,$ if $s<k+1,$ where $A$ is microlocal cutoff of order zero such that $Au\in H^{s}.$
$\Box$


Now suppose $s\geq m>1/2-\text{Im }\lambda$, $-k<m< s\leq k+1,$ $u\in H^{m}(X_{-\delta_0})$ and $(P_{\lambda}-iQ_{\lambda})u\in H^{s-1}(X_{-\delta_0})$. Here we consider a smaller interval, namely $-k<s<k+1$ instead of $s<k+1$ because of a duality condition that will be explained latter. We have shown above that $\text{WF}^{s}(u)$ (actually, more than that, $\text{WF}^{s+1}(u)$) is disjoint from the elliptic set of $P_{\lambda}-iQ_{\lambda}$. Because of the properties of the support of $Q_{\lambda}$, $\Lambda_{\pm}\subset\text{ell}(Q_{\lambda})$ is disjoint from $\text{WF}^{s}(u)$ and since $\text{WF}^{s}(u)$ is a closed set, there is a neighborhood of $\Lambda_{\pm}$ disjoint from $\text{WF}^{s}(u).$ Then, it follows from \eqref{ellq}, that $\text{WF}^{s}(u)\cap\Sigma_{\pm}=\emptyset.$ Then by Proposition \ref{propsing}, we have that $u\in H^{s}(X_{-\delta_0})$. 

It then makes sense to consider the inclusion map
$$i_{s}:E_{s}=\{u\in H^{m}(X_{-\delta_0}):(P_{\lambda}-iQ_{\lambda})u\in H^{s-1}(X_{-\delta_0})\}\rightarrow H^{s}(X_{-\delta_0}).$$
As it is showed in \cite{V2}, one can show that $E_{s}$ is a complete Banach space with the norm $\displaystyle\|u\|_{E_s}^2=\|u\|_{H^m}^2+\|(P_{\lambda}-iQ_{\lambda})u\|_{H^{s-1}}.$

Also the graph of the inclusion map is closed. Thus, by the closed graph theorem, the map $i_{s}:E_{s}\rightarrow H^{s}(X_{-\delta_0})$ is continuous, that is
\begin{equation}\label{estinclu}
\|u\|_{H^{s}}\leq C\left( \|(P_{\lambda}-iQ_{\lambda})u\|_{H^{s-1}}+\|u\|_{H^m}\right),\quad u\in E_s.
\end{equation}

The estimate \eqref{estinclu} implies that if $u\in\text{Ker}(P_{\lambda}-iQ_{\lambda})\cap H^s(X_{-\delta_0})$, then $\|u\|_{H^s}\leq \|u\|_{H^m}$. Since $X_{-\delta_0}$ is a pre-compact space without boundary and $H^m(X_{-\delta_0})\subset H^s(X_{-\delta_0})$, the inclusion map $H^{s}(X_{-\delta_0})\rightarrow H^m(X_{-\delta_0})$ is compact by the Rellich-Kondrachov theorem. Then the unit ball of $\text{Ker}(P_{\lambda}-iQ_{\lambda})$ is compact, which implies that $\text{Ker}(P_{\lambda}-iQ_{\lambda})$ has finite dimension. Moreover, it follows from Proposition \ref{propsing} that $\displaystyle \text{Ker}(P_{\lambda}-iQ_{\lambda})\subset \cap_{0<r<k+1}H^{r}(X_{-\delta_0})$, and thus this space is independent of the choice of $s$, $s<k+1.$ 

Now, we consider the adjoint operator $P_{\lambda}^{*}+iQ_{\lambda}^{*}$. We have that if $s'<1/2+\text{Im }\sigma,$ $-k<s'<k+1,$ $u\in H^{m}(X_{-\delta_0})$ and $(P_{\lambda}^{*}+iQ_{\lambda}^{*})u\in H^{s'-1}(X_{-\delta_0})$, then by similar arguments as the ones used above, using now Proposition \ref{propsinga}, $u\in H^{s'}(X_{-\delta_0})$. It then makes sense to consider the map
$$j_{s'}:B_{s'}=\{u\in H^{m}(X_{-\delta_0}):(P_{\lambda}^{*}+iQ_{\lambda}^{*})u\in H^{s'-1}(X_{-\delta_0})\}\rightarrow H^{s'}(X_{-\delta_0}).$$

Analogously as above, we can conclude that 

\begin{equation}\label{estincla}
\|u\|_{H^{s'}}\leq C(\|(P_{\lambda}^{*}+iQ_{\lambda}^{*})u\|_{H^{s'-1}}+\|u\|_{H^{m}}),\quad u\in B_{s'},
\end{equation}
which implies that $\text{Ker}(P_{\lambda}^{*}+iQ_{\lambda}^{*})$ in $H^{s'}(X_{-\delta_0})$ has finite dimension. 

We now follow similar arguments as the ones given in \cite{V2}, which, as pointed out in \cite{V2}, come from \cite{Ho} (proof of Theorem 26.1.7), and show that estimate \eqref{estincla} gives the $H^{s}-$solvability of 
$$(P_{\lambda}-iQ_{\lambda})u=f,\quad s>1/2-\text{Im }\lambda, -k<s<k+1,$$
for $f$ in the annihilator of $\text{Ker}(P_{\lambda}^{*}+iQ_{\lambda}^{*})$, that is $f\in \text{Ran}_{H^{s}}(P_{\lambda}-iQ_{\lambda})$. This shows that the range of $P_{\lambda}-iQ_{\lambda}$ in $H^{s}(X_{-\delta_0})$ is closed. This in turn with the fact that its kernel and cokernel in $H^{s}(X_{-\delta_0})$ are finite-dimensional implies that $P_{\lambda}-iQ_{\lambda}$ is an analytic family of Fredholm operators for $s<1/2-\text{Im }\lambda$, $-k<s<k+1.$ We remark here that in Theorem 26.1.7 of \cite{Ho} it is assumed that the manifold $X_{-\delta_0}$ is non-trapping (that is, no complete bicharacteristic curve is contained in $X_{-\delta_0}$), but, as remarked in \cite{Ho}, this condition is merely sufficient to prove that the kernel of $P_{\lambda}^{*}+iQ_{\lambda}^{*}$ is finite-dimensional, and this hypothesis is no way necessary for the conclusion to be valid. We actually showed above that both $\text{Ker}(P_{\lambda}^{*}+iQ_{\lambda}^{*})$ and $\text{Ker}(P_{\lambda}-iQ_{\lambda})$ have finite dimension without assuming this condition.

Now let us prove that the range $P_{\lambda}-iQ_{\lambda}$ in $H^{s}(X_{-\delta_0})$ is closed. The dual of $H^{s}$ for $s>1/2-\text{Im }\lambda,$ $-k<0\leq s<k+1$, is $H^{-s}(X_{-\delta_0})=H^{s'-1}(X_{-\delta_0})$ with $s'=1-s,$ so $s'<1/2+\text{Im }\lambda$, $-k<s'<k+1.$ First, let $V$ be a subspace of $H^{s'}$ in the orthogonal complement of $\text{Ker}(P_{\lambda}^{*}+iQ_{\lambda}^{*})$, that is $V$ is closed, $V\cap\text{Ker}(P_{\lambda}^{*}+iQ_{\lambda}^{*})=\{0\}$ and $V+\text{Ker}(P_{\lambda}^{*}+iQ_{\lambda}^{*})=H^{s'}(X_{-\delta_0}).$ Then for $v\in V\cap B_{s'}$, we have that
\begin{equation}\label{veagasl}
\|v\|_{H^{s'}}\leq C'\|(P_{\lambda}^{*}+iQ_{\lambda}^{*})v\|_{H^{s'-1}},\quad v\in V\cap B_{s'}. 
\end{equation}

Indeed, suppose that the inequality \eqref{veagasl} is not valid for some element of $V\cap B_{s'}$. Then one can construct a sequence $\{u_{j}\}\in V\cap B_{s'}$ such that $\|u_{j}\|_{H^{s'}}=1$ and $\|(P_{\lambda}^{*}+iQ_{\lambda}^{*})u_j\|_{H^{s'-1}}\rightarrow 0$. By the weak compactness of the unit ball in $H^{s'}(X_{-\delta_0})$, there exists a subsequence $\{u_{jk}\}$ which converges weakly to some $u\in H^{s'}(X_{-\delta_0})$. Since $V$ is closed, $u\in V$, thus $(P_{\lambda}^{*}+iQ_{\lambda}^{*})u_{jk}\rightarrow (P_{\lambda}+iQ_{\lambda}^{*})u$ weakly in $H^{s'-2}(X_{-\delta_0})$. Hence $(P_{\lambda}^{*}+iQ_{\lambda}^{*})u=0$, since $H^{s'-2}(X_{-\delta_0})\subset H^{s'-1}(X_{-\delta_0}).$ This implies that $u\in V\cap \text{Ker}(P_{\lambda}^{*}+iQ_{\lambda}^{*})=\{0\}.$ However, since the inclusion $H^{s'}(X_{-\delta_0})\rightarrow H^{-N}(X_{-\delta_0})$ is compact, $u_{jk}$ converges strongly to $u$ in $H^{-N}(X_{-\delta_0})$, so by \eqref{estincla} $u_{jk}\rightarrow u$ in $H^{s'}(X_{-\delta_0})$. This implies in particular that $\|u\|_{H^{s'}}=1,$ which contradicts $u=0$. 

Now, for $s'=1-s,$ for $f$ in the annihilator of $\text{Ker}(P_{\lambda}^{*}+iQ_{\lambda}^{*})\subset H^{s'}(X_{-\delta_0})$, and for $v\in V\cap B_{s'},$ \eqref{veagasl} implies that
$$|\left\langle f,v\right\rangle|\leq \|f\|_{H^{s-1}}\|v\|_{H^{s'}}\leq C'\|f\|_{H^{s-1}}\|(P_{\lambda}^{*}+iQ_{\lambda}^{*})v\|_{H^{s'-1}},$$
which implies that for any $w\in\text{Ker}(P_{\lambda}^{*}+iQ_{\lambda}^{*})$
$$|\left\langle f,v+w\right\rangle|=|\left\langle f,v\right\rangle|\leq  C'\|f\|_{H^{s-1}}\|(P_{\lambda}^{*}+iQ_{\lambda}^{*})v\|_{H^{s'-1}}$$
$$\text{then }|\left\langle f,v+w\right\rangle|\leq C'\|f\|_{H^{s-1}}\|(P_{\lambda}^{*}+iQ_{\lambda}^{*})(v+w)\|_{H^{s'-1}}.$$

Therefore the inequality \eqref{veagasl} holds for all $v\in B_{s'}\subset H^{s'}.$ So the conjugated linear map
$$(P_{\lambda}^{*}+iQ_{\lambda}^{*})v\mapsto \left\langle f,v\right\rangle,\quad v\in B_{s'}$$
is well-defined and continuous from $\text{Ran}_{B_{s'}}(P_{\lambda}^{*}+iQ_{\lambda}^{*})\subset H^{s'-1}$ to $\mathbb{C}$. It then follows from the Hahn-Banach theorem that this map can be extended to a conjugated linear functional $\omega$ on $H^{s'-1}$, and by Riesz representation theorem there exists $u\in H^{s}$ such that $\left\langle u,\varphi\right\rangle=\omega(\varphi)$ for any $\phi\in H^{s'-1}.$ In particular, for $\varphi=(P_{\lambda}^{*}+iQ_{\lambda}^{*})\psi,$ $\psi\in C^{\infty}(X)\subset B_{s'}$, 
$$\left\langle u,(P_{\lambda}^{*}+iQ_{\lambda}^{*})\psi\right\rangle=\left\langle f,\psi\right\rangle.$$

Since $\psi$ is arbitrary and $C^{\infty}(X)$ is dense in $B_{s'}$, $(P_{\lambda}-iQ_{\lambda})u=f,$ as we wanted to show. This proves Theorem \ref{meroextention}.

\section{Semiclassical Estimates}

The analysis of the semiclassical problem is analogous to the one in the classical setting. Again, as before, by elliptic regularity, for $-k<\delta<s\leq k,$ 
$$P_{h,z}-iQ_{h,z}\text{ elliptic at }\nu,\nu\not\in \text{WF}_{\hbar}^{s-2,0}((P_{h,z}-iQ_{h,z})u)\cup\text{WF}_{\hbar}^{\delta,-1}(u)\Rightarrow \nu\not\in \text{WF}_{\hbar}^{s,0}(u).$$
Then if $(P_{h,z}-iQ_{h,z})u\in H^{s-2}$ and $p_{\hbar,z}-iq_{\hbar,z}$ is elliptic at $\nu$ implies that $\nu\not\in \text{WF}_{\hbar}^{s,0}(u).$ Also, similarly as before, if $-k<s\leq k$, then if $\nu\not\in \text{WF}_{\hbar}^{s,-1}(u)$ and $\displaystyle \tilde{\gamma}(\nu)\cap\text{WF}_{\hbar}^{s-1,0}((P_{h,z}-iQ_{h,z})u)=\emptyset,$ with $\tilde{\gamma}(\nu)$ a forward (+) or backward (-) bicharacteristic from $\nu$, since $P_{h,z}-iQ_{h,z}$ is elliptic and $\pm q_{\hbar,z}\geq 0$ on a neighborhood of $\tilde{\gamma}_{\pm}(\nu),$ then $\displaystyle \tilde{\gamma}(\nu)\cap \text{WF}_{\hbar}^{s,-1}(u)=\emptyset.$ Again, for $P_{h,z}^{*}+iQ_{h,z}^{*}$ the directions are reversed. 

We first prove the following semiclassical regularity result. 

\begin{prop}\label{psemi}
Suppose $s>m>1/2-\text{Im }\lambda$, $\displaystyle m< s< k+1$, and $\text{WF}_{\hbar}^{m,-1}(u)\cap L_{\pm}=\emptyset.$ Then
$$L_{\pm}\cap\text{WF}^{s-1,0}(P_{\hbar,z}u)=\emptyset\Rightarrow  L_{\pm}\cap\text{WF}^{s,-1}(u)=\emptyset.$$
\end{prop}
\textbf{Proof.} The proof is analogous to the one given in Proposition \ref{propsing}, but here we need to choose a $\phi\in C_{0}^{\infty}(\mathbb{R})$ that also localizes in $\tilde{\rho}$. Let $\phi\in C_{0}^{\infty}(\mathbb{R})$ be such that $\phi$ is $1$ in a neighborhood of $0$, $\phi'(t)\leq 0$ for all $t\geq0$ and $\phi$ supported sufficiently close to $0$ so that 
\begin{equation}\label{hamisemiine}
\tilde{\rho}^2+\rho_0\in\text{supp }\phi'\Rightarrow (\sgn \xi)\tilde{\rho}H_{\text{Re }p_{\hbar,z}}(\tilde{\rho}^2+\rho_0)>0.
\end{equation}
Such $\phi$ exists since  $\displaystyle(\sgn \xi)\tilde{\rho}H_{\text{Re }p_{\hbar,z}}(\tilde{\rho}^2+\rho_0)\geq 8(\tilde{\rho}^2+\rho_0)-\mathcal{O}((\tilde{\rho}^2+\rho_0)^{3/2}).$ Indeed,
$$(\sgn \xi)\tilde{\rho}H_{\text{Re }p_{\hbar,z}}(\rho_{0}+\tilde{\rho}^2)=\left(8+\frac{4\tilde{\rho}\text{Re }z}{(1+\mu)^2}+4\tilde{\rho}^{2}\partial_{\mu}|\eta|^{2}\right)(\rho_0+\tilde{\rho}^2+\tilde{\rho}^4)$$
Then let $c$ be given by 
$$ c=\phi(\rho_0+\tilde{\rho}^2)\tilde{\rho}^{-s+1/2},\quad c_{\epsilon}=c(1+\epsilon\tilde{\rho}^{-1})^{-\delta}.$$

If we then proceed as in the proof of Proposition \ref{propsing}, we have that 
$$P_{\hbar,z}^{*}C_{\epsilon}^{*}C_{\epsilon}-C_{\epsilon}^{*}C_{\epsilon}P_{\hbar,z}=(P_{\hbar,z}^{*}-P_{\hbar,z})C_{\epsilon}^{*}C_{\epsilon}+[P_{\hbar,z},C_{\epsilon}^{*}C_{\epsilon}].$$

Now write $\displaystyle P_{\hbar,z}=P_{\hbar,z,Re}+iP_{\hbar,z,Im},$ where 
$$ P_{\hbar,z,Re}=\frac{P_{\hbar,z}+P_{\hbar,z}^{*}}{2},\quad P_{\hbar,z,Im}=\frac{P_{\hbar,z}-P_{\hbar,z}^{*}}{2i},$$
and the principal symbols of $\displaystyle P_{\hbar,z,Re}$ and $P_{\hbar,z,Im}$ are given by (\ref{realsymb}) and (\ref{imasemi}), respectively. We have that
\begin{equation}\label{semicommu}
(\sgn \xi)\frac{i}{h}(P_{\hbar,z}^{*}C_{\epsilon}^{*}C_{\epsilon}-C_{\epsilon}^{*}C_{\epsilon}P_{\hbar,z})=(\sgn \xi)\left(\frac{i}{h}[P_{\hbar,z,Re},C_{\epsilon}^{*}C_{\epsilon}]+\frac{1}{h}(P_{\hbar,z,Im}C_{\epsilon}^{*}C_{\epsilon}+C_{\epsilon}^{*}C_{\epsilon}P_{\hbar,z,Im})\right).
\end{equation}
Then
$$\left\langle ih^{-1}(P_{\hbar,z}^{*}C_{\epsilon}^{*}C_{\epsilon}-C_{\epsilon}^{*}C_{\epsilon}P_{\hbar,z})u,u\right\rangle=\left\langle ih^{-1}[P_{\hbar,z,Re},C_{\epsilon}^{*}C_{\epsilon}]u,u\right\rangle+2h^{-1}\text{Re}\left(\left\langle P_{\hbar,z,Im}u,C_{\epsilon}^{*}C_{\epsilon}u\right\rangle\right).$$
Now, when $\mu=0,\eta=0$, the semiclassical principal symbol of $\displaystyle h^{-1}i[P_{\hbar,z,Re},C_{\epsilon}^*C_{\epsilon}]$ is given by
\begin{equation}\label{princesemireal}
\begin{gathered}
(8(\sgn \xi)+4\tilde{\rho}\text{Re }z)\left(\frac{1}{2}-s+\frac{\delta\epsilon}{\tilde{\rho}+\epsilon}\right)\phi^{2}\tilde{\rho}^{-2s}(1+\epsilon\tilde{\rho}^{-1})^{-2\delta}\\+2(\sgn \xi)\phi'\phi\tilde{\rho}\left(H_{\text{Re }p_{\hbar,z}}(\rho_{0}+\tilde{\rho}^2)\right)\tilde{\rho}^{-2s}(1+\epsilon\tilde{\rho}^{-1})^{-2\delta}.
\end{gathered}
\end{equation}

We can write 
$$\text{Re}\left(\left\langle P_{\hbar,z,Im}u,C_{\epsilon}^{*}C_{\epsilon}u\right\rangle\right)=\text{Re}\left(\left\langle C_{\epsilon}P_{\hbar,z,Im}u,C_{\epsilon}u\right\rangle\right)=\left\langle P_{\hbar,z,Im}C_{\epsilon}u,C_{\epsilon}u\right\rangle-\text{Re}\left(\left\langle [P_{\hbar,z,Im},C_{\epsilon}]u,C_{\epsilon}u\right\rangle\right).$$

When $\mu=0,\eta=0$, the semiclassical principal symbol of $\displaystyle C_{\epsilon}^*P_{\hbar,z,Im}C_{\epsilon}$ is 
$$c_{\epsilon}^2\text{Im }p_{\hbar,z}=-2(\text{Im }z)[2(\sgn \xi)\tilde{\rho}^{-1}+\text{Re }z]\phi^2\tilde{\rho}^{-2s+1}(1+\epsilon\tilde{\rho}^{-1})^{-2\delta}.$$

Finally, note that the semiclassical principal symbol of $h^{-1}C_{\epsilon}^*[P_{\hbar,z,Im},C_{\epsilon}]$ is $-ic_{\epsilon}H_{\text{Im }p_{\hbar,z}}c_{\epsilon},$ which is in $S^{2s-1}(k+1-2s)$ and is purely imaginary. Therefore the symbol of $\text{Re }(h^{-1}C_{\epsilon}^*[P_{\hbar,z,Im},C_{\epsilon}])$ is in $hS^{2s-2}(k+1-2s).$ We then have that, if $\mu=0,\eta=0$, using $h^{-1}z=\lambda,$ the semiclassical principal symbol of (\ref{semicommu}) is given by

\begin{equation}\label{princesemi}
\begin{gathered}
(8+4(\sgn \xi)\tilde{\rho}\text{Re }z)\left(-\text{Im }\lambda+\frac{1}{2}-s+\frac{\delta\epsilon}{\tilde{\rho}+\epsilon}\right)\phi^{2}\tilde{\rho}^{-2s}(1+\epsilon\tilde{\rho}^{-1})^{-2\delta}\\+2(\sgn \xi)\phi'\phi\tilde{\rho}\left(H_{\text{Re }p_{\hbar,z}}(\rho_{0}+\tilde{\rho}^2)\right)\tilde{\rho}^{-2s}(1+\epsilon\tilde{\rho}^{-1})^{-2\delta}.
\end{gathered}
\end{equation}
Note that when $\tilde{\rho}$ is small, (\ref{princesemi}) is less than or equal to 
$$-8\left(\text{Im }\lambda-\frac{1}{2}+s-\delta\right)\phi^{2}\tilde{\rho}^{-2s}(1+\epsilon\tilde{\rho}^{-1})^{-2\delta}$$

Indeed, this follows if we have $\displaystyle s>1/2-\text{Im }\lambda,$ $\phi'\phi<0$ and (\ref{hamisemiine}). We now proceed exactly as in Proposition \ref{propsing}.
$\Box$

Now, for $P_{\lambda}^{*}$ the proof is an easy adaptation of Proposition \ref{psemi}, following a similar proof as the one given in Proposition \ref{propsinga}.
\begin{prop}
For $s<1/2+\text{Im }\lambda$, $\displaystyle  s\leq k+1$, and $O$ a neighborhood of $L_{\pm}$
$$\text{WF}^{s,-1}_{\hbar}(u)\cap(O\setminus L_{\pm})=\emptyset,\text{ WF}_{\hbar}^{s-1,0}(P_{\lambda}^{*}u)\cap L_{\pm}=\emptyset\Rightarrow
\text{WF}_{\hbar}^{s,-1}(u)\cap L_{\pm}=\emptyset.$$
\end{prop}

We now assume that $p_{\hbar,z}$ is semiclassicaly non-trapping. Suppose that $s\geq m>1/2-\text{Im }\lambda$, $-k<s\leq k+1,$ $h^{N}u\in H_{\hbar}^m$ and $(P_{h,z}-iQ_{h,z})u\in H_{\hbar}^{s-1}.$ Then analogously as in the classical setting, using the properties of the complex absorption and the semiclassical regularity result (Proposition \ref{psemi}), one can conclude that $hu\in H_{\hbar}^{s}.$ We also have a similar result, as we did in the previous section, for the adjoint $P_{h,z}^{*}+iQ_{h,z}^{*}$.

Differently from the classical problem, it is easier to prove that the operator $P_{h,z}-iQ_{h,z}:E_{s}\cap H^{s}\rightarrow H^{s}$ is Fredholm. The following lemma is proved in \cite{V2}.

\begin{lemma}
$\displaystyle \text{Ker}(P_{h,z}-iQ_{h,z})=\{0\}=\text{Ker}(P_{h,z}^{*}+iQ_{h,z}^{*}),$ for $h$ sufficiently small, as operators acting on $H_{\hbar}^{s}$.
\end{lemma}


Analogously as in the classical case shown previously, we can deduce that the range of $P_{h,z}-iQ_{h,z}$ is closed. This shows that $P_{h,z}-iQ_{h,z}:E_{s}\cap H^{s}\rightarrow H^{s}$ is Fredholm, which implies the existence of a $h_0$ such that $P_{h,z}-iQ_{h,z}$ is invertible for $h<h_0,$ $s\geq m>1/2-\text{Im }\lambda$, $-k<s\leq k+1.$

Now, to obtain uniform estimates of $(P_{h,z}-iQ_{h,z})^{-1}$ as $h\rightarrow0,$ we normalize the problem by introducing weights in order to make the function spaces independent of $h$. Then for $r\in\mathbb{R},$ let $D_{h}^{r}\in \Psi_{h}^{r}$ be an elliptic and invertible operator. We let
$$P_{h,z}^{s}-iQ_{h,z}^{s}=D_{h}^{s-1}(P_{h,z}-iQ_{h,z})D_{h}^{s}.$$
Then, for a fixed $h_0$ and $z_0$, let $\displaystyle\mathcal{E}=\{u\in L^2:(P_{h_0,z_0}^{s}-iQ_{h_0,z_0}^{s})u\in L^{2}\}.$ Let $i:\mathcal{E}\rightarrow L^2$ be the inclusion map. We have that the map
$$i\circ h(P_{h,z}-iQ_{h,z}):\mathcal{E}\rightarrow L^2$$
is continuous for any $h<h_0.$ 

Now let $f\in L^2$, and $v_h=h(P_{h,z}^s-iQ_{h,z}^s)^{-1}f.$ We will show that $v_h$ is uniformly bounded in $h$. Indeed, suppose first that $hv_h$ is not bounded. Then there exists a sequence $\{h_j\}$ such that $h_j\|v_{h_j}\|_{L^2}\geq 1. $ Then let 
$$w_h=\frac{v_h}{h^2\|v_h\|_{L^2}},\quad h\in \{h_j\}.$$

We have that $\|w_h\|_{L^2}=h^{-2},$ and since $h_j\|v_{h_j}\|_{L^2}\geq 1$
$$(P_{h,z}^s-iQ_{h,z}^s)w_h=\frac{f}{h\|v_h\|_{L^2}}\Rightarrow \|(P_{h,z}^s-iQ_{h,z}^s)w_h\|_{L^2}=\frac{\|f\|_{L^2}}{h\|v_h\|_{L^2}}\leq \|f\|_{L^2}.$$

Then $h^2 w_h$ is uniformly bounded in $H_{\hbar}^{0}$ and $(P_{h,z}^{s}-iQ_{h,z}^{s})w_h\in H^{0}_{\hbar}=L^2_{\hbar},$ so by what we discussed above $h w_h\in L^2.$ That is, 
$\displaystyle \frac{v_h}{h\|v_h\|_{L^2}}$, which has $L^2-$norm $h^{-1}$, is bounded, an absurd. Hence $hv_h\in L^2.$

Now let $u_h=h^{-1}v_{h}$, so $h^{2}u_h\in L^2$, and $(P_{h,z}^{s}-iQ_{h,z}^{s})u_h=f\in L^2.$ Then again by the same arguments as before, $hu_h\in L^2.$ Therefore, by the uniform boundedness principle, $j\circ h(P_{h,z}^{s}-iQ_{h,z}^{s})^{-1}$ is equicontinuous. This implies that
$$\|(P_{h,z}-iQ_{h,z})^{-1}f\|_{H_{\hbar}^s}\leq Ch^{-1}\|f\|_{H_{\hbar}^s}.$$

As discussed at the end of the section where we introduced the complex absorption, the assumption of semiclassical non-trapping is valid, for instance, in any region where the operator $P_{\lambda}-iQ_{\lambda}$ is elliptic. In particular, as we saw before, $p_{\hbar,z}-iq_{\hbar,z}$ is elliptic in conic sectors $\delta|\text{Re }\lambda|< \text{Im }\lambda<\delta_0|\text{Re }\lambda|,$ for some $\delta_0>0,$ $0<\delta<\delta_0.$ So, using $\lambda=h^{-1}z,$ $h=|\lambda|^{-1}$, we have proved above the following theorem. 

\begin{thm}\label{thmpq}
Let $s>1/2-\text{Im }\lambda,$ $-k<s\leq k+1$. Then $P_{\lambda}-iQ_{\lambda}:E_{s}\rightarrow H^{s-1}$ admits a meromorphic inverse $(P_{\lambda}-iQ_{\lambda})^{-1}:H^{s-1}\rightarrow E_{s}$. Moreover, if the metric $g$ is non-trapping, the following estimate holds in $-C<\text{Im }\lambda<\delta_0|\text{Re }\lambda|,$ $|\text{Re }\lambda|>C$ 
\begin{equation}\label{nontrapest}
\|(P_{\lambda}-iQ_{\lambda})^{-1}f\|_{H_{|\lambda|^{-1}}^s}\leq C'|\lambda|^{-1}\|f\|_{H_{|\lambda|^{-1}}^{s-1}},
\end{equation}
\end{thm}

Then the following theorem follows directly from \eqref{resolpq} and Theorems \ref{meroextention} and \ref{thmpq}. 

\begin{thm}
Let $(X,g)$ be an asymptotically hyperbolic manifold, $x$ a boundary defining function on $X$ and $k\in \mathbb{N}\cup\{\infty\}.$ Then if $g$ is even modulo $O(x^{2k+1})$ the modified resolvent 
$$R(\lambda)=\left(\Delta_g-\frac{n^2}{4}-\lambda^2\right)^{-1}:\mathcal{\dot{S}}\rightarrow \mathcal{\dot{S}}'$$
defined for $\text{Im }\lambda>>0$, extends to a finite-meromorphic family from $\text{Im }\lambda>>0$ to $\text{Im }\lambda>-1/2-k.$ Moreover, if the metric $g$ is non-trapping, the following estimate holds in $1/2-s<-C<\text{Im }\lambda<\delta_0|\text{Re }\lambda|,$ $|\text{Re }\lambda|>C$ 
\begin{equation}\label{nontrapest2}
\|x^{-n/2}e^{i\lambda\phi}R(\lambda)f\|_{H_{|\lambda|^{-1}}^s}\leq C'|\lambda|^{-1}\|x^{-(n+4)/2}e^{i\lambda\phi}f\|_{H_{|\lambda|^{-1}}^{s-1}},
\end{equation}
where $e^{\phi}=\mu^{1/2}(1+\mu)^{-1/4}$ near $\{\mu=0\}.$
\end{thm}

\end{document}